\documentclass[12pt]{article}
\usepackage{amsmath,amsfonts,amssymb,amsthm,amscd}

\catcode`\@=11
\renewcommand\section{\@startsection{section}{1}{\z@}%
                                      {-3.25ex\@plus -1ex \@minus -.2ex}%
                                      {1.5ex \@plus .2ex}%
                                      {\normalfont\large\bfseries}}
\catcode`\@=\active
\catcode`\@=11
\renewcommand\subsection{\@startsection{subsection}{2}{\z@}%
                                      {-3.25ex\@plus -1ex \@minus -.2ex}%
                                      {1.5ex \@plus .2ex}%
                                      {\normalfont\large\bfseries}}
\catcode`\@=\active

\newcommand{\w}{\omega}
\newcommand{\1}{{\bf 1}}
\newcommand{\vactw}{1_{\rm tw}}
\newcommand{\id}{{\rm id}}
\newcommand{\Hom}{{\rm Hom\,}}
\DeclareMathOperator{\Ker}{Ker\,}
\DeclareMathOperator{\Aut}{Aut\,}
\newcommand{\End}{{\rm End\,}}
\newcommand{\wt}[1]{{\rm wt}(#1)}

\newcommand{\Ext}[1]{{\rm Ext}^1_{#1}\,}
\newcommand{\image}{{\rm Im\,}}
\newcommand\Z{\mathbb{Z}}
\newcommand\Zpos{\Z_{\geq0}}
\newcommand\Q{\mathbb{Q}}
\newcommand\C{\mathbb{C}}

\newcommand\h{\mathfrak{h}}
\newcommand\g{\mathfrak{g}}

\newcommand\hh{\hat{\h}}

\newcommand{\NO}{\,{\raise0.25em\hbox{$\mathop{\hphantom {\cdot}}\limits^{_{\circ}}_{^{\circ}}$}}\,}
\newcommand{\Free}[1]{M(1)^{#1}}
\newcommand{\Fremo}[1]{M(1,#1)}
\newcommand{\charge}[1]{V_{L}^{#1}}

\newcommand{\charlam}[1]{V_{#1+L}}
\newcommand{\Fretw}[1]{M(1)(\theta)^{#1}}

\newcommand{\Tl}[2]{\widetilde{#1}(#2)}

\newcommand\haru[2]{{\rm span}\{\,#1\,|\,#2\,\}}
\newcommand{\fusion}[3]{{\binom{#3}{#1\;#2}}}

\newtheorem{theorem}{Theorem}[section]
\newtheorem{proposition}[theorem]{Proposition}
\newtheorem{lemma}[theorem]{Lemma}
\newtheorem{corollary}[theorem]{Corollary}
\theoremstyle{definition}
\newtheorem{definition}[theorem]{Definition}

\theoremstyle{remark}
\newtheorem{remark}[theorem]{Remark}

\numberwithin{equation}{section}

\makeatletter
\@addtoreset{equation}{section}

\makeatother

\begin{document}
\begin{center}
\textbf{\Large Rationality of the vertex operator algebra $\charge{+}$\\
 for a positive definite even lattice $L$}
\end{center}

\vskip 2ex
\begin{center}
Toshiyuki Abe\footnote{Supported by JSPS Research Fellowships for Young 
Scientists.}\\

\vskip 2ex
Graduate School of Mathematical Sciences , University of Tokyo\\
3-8-1 Komaba, Meguroku, Tokyo 153-8914, Japan\\
{\rm e-mail: abe@ms.u-tokyo.ac.jp}
\end{center}

\vskip 3ex
\noindent
\begin{small}
\textbf{Abstract:}
The lattice vertex operator algebra $V_L$ associated to a positive definite even lattice $L$ has an automorphism of order $2$ lifted from $-1$-isometry of $L$. 
We prove that for the fixed point vertex operator algebra $V_L^+$, any $\Zpos$-graded weak module is completely reducible.    
\end{small}

\section{Introduction}
Let $\charge{}$ be a lattice vertex operator algebra associated with a positive definite even lattice $L$. 
The vertex operator algebra $\charge{}$ has an automorphism $\theta$ of order $2$ lifted from the $-1$-isometry of $L$, and the fixed point set $\charge{+}$ is naturally a vertex operator algebra. 
In the paper, we prove the semisimplicity of arbitrary $\Zpos$-graded $\charge{+}$-module for any rank one positive definite even lattice. 
This fact has been proved in \cite{A3} under the restricted condition on the lattice $L$.  

A vertex operator algebra $V$ is called \textit{rational} if any $\Zpos$-graded $V$-module is completely reducible. 
Many examples of rational vertex operator algebras are known; vertex operator algebras associated with even lattices (\cite{FLM}, \cite{D1}), vertex operator algebras of irreducible vacuum representations for affine Kac-Moody algebras with positive levels (\cite{FZ} and \cite{L3}), vertex operator algebras of irreducible highest weight modules for the Virasoro algebra of highest weight zero and minimal central charges (\cite{FZ} and \cite{W1}). 
These rational vertex operator algebras satisfy at the same time $C_2$-condition, which is the property that the subspace $C_2(V):=\{a(-2)b\,|\,a,b\in V\}$ is of finite codimension in $V$. 
In \cite{ABD}, it is proved that a rational vertex operator algebra satisfying $C_2$-condition is regular, i.e., any weak module is semisimple. 
In particular, the rational vertex operator algebras above are all regular.
The $C_2$-condition plays a central role in the conformal field theory on vertex operator algebras (cf. \cite{NT} and \cite{H1}--\cite{H4}). 
But one needs further conditions to develop the structure of conformal blocks in more detail.
Rationality together with $C_2$-condition enjoy all necessary conditions in their developments and therefore regular vertex operator algebras give mathematical models in conformal field theories. 

For a vertex operator algebra $V$ and a finite group $G$ of automorphisms of $V$, the fixed point set $V^{G}=\{\,a\in V\,|\,ga=a\hbox{ for }g\in G\,\}$ has a natural vertex operator algebra structure. 
Structure of $V$-modules and twisted modules as $V^G$-module have been studied in \cite{DM}, \cite{HMT}, \cite{DY} and \cite{MT}.
However, it is still a conjecture that any irreducible $V^{G}$-module appears in irreducible $V$-modules or twisted modules as a $V^{G}$-submodule.  
This conjecture was proved for the lattice vertex operator algebra with the group $\langle\theta\rangle$ in \cite{DN2} for the rank one case and in \cite{AD} for higher rank lattices. 
There are another conjectures: (I) Is $V^G$ rational if $V$ is rational? (II) Does $V^G$ satisfy the $C_2$-condition if $V$ satisfies the condition?
Conjecture (II) for $\charge{+}:=\charge{\langle\theta\rangle}$ is solved by \cite{Gail} in the case rank one and by \cite{ABD} in the case general rank. 
Our result solves conjecture (I) affirmatively for $\charge{}$ with the automorphism group $G=\langle\theta\rangle$.         

Conjecture (I) for $\charge{}$ with arbitrary finite automorphism group $G$ and  a rank one lattice $L$ is also related to the classification problem of the rational vertex operator algebras with $c=1$ (\cite{Gin}, \cite{DG} and \cite{DGR}). 
The partition functions with modular invariance of rational conformal field theory with $c=1$ are considered in \cite{Gin} and classified by Kiritsis in \cite{Kiritsis}.
They are given by the characters of the fixed point vertex operator algebras $\charge{G}$ of $\charge{}$ with rank one even lattice $L$ by finite automorphism groups $G$.
So we expect that the characters are of rational and $C_2$-cofinite vertex operator algebras motivated form \cite{Zh1}.   
The result of the present paper proves this for $c=1$ rational conformal field theories associated with $G=\langle\theta\rangle$. 
The remaining cases occur only in the case $L=\Z\sqrt{2}h$ with $(h,h)=1$ (see \cite{DGR}). 
But the classification of irreducible $\charge{G}$-modules has not been completed yet, and it is still open that $\charge{G}$ is rational and satisfies the $C_2$-condition.  
 
We explain the way to prove the rationality of $\charge{+}$.
Since the proof of the rationality reduces to the case of rank one by using Theorem 6.11 in \cite{Miya3} (see Section \ref{Generalrank}), we mainly consider the case of rank one. 

Let $L$ be a rank one positive definite even lattice. 
The main claim which we prove is that the extension $\Ext{\charge{+}}(M^2,M^1)=0$ for any irreducible $\charge{+}$-modules $M^{1},\,M^2$.
This fact with $C_2$-finiteness of $\charge{+}$ shows that any $\Zpos$-graded weak $\charge{+}$-module is completely reducible. 
The proof of the triviality of extensions for a pair $(M^1,M^2)$ is divided into  the four cases.
First case is the case that the difference of the lowest weights of $M^1$ and $M^2$ is not an integer, where the triviality is clear. 
Second one is the case the lowest weights are same. 
In this case we use the semisimplicity of Zhu's algebra associated to $\charge{+}$ proved in \cite{DN2}.
Third one is the case $(M^1,M^2)=(\charge{\pm},\charge{\mp})$. 
In the case the triviality of $\Ext{\charge{+}}(M^2,M^1)=0$ is proved in \cite{A3} essentially. 
Fourth one is the remaining cases. 
The triviality of extensions in the fourth case follows from the representation theory for $\Free{+}$ (see \cite{DG} and \cite{DN1}).
We discuss the splitness of a short exact sequence $0\rightarrow N^1\rightarrow N\rightarrow N^2\rightarrow0$ for irreducible $\Free{+}$-modules $N^1,\,N^2$ and a weak $\Free{+}$-module $N$. 
Since $\Free{+}$ is not rational, every short exact sequences are not splitting.
But as in Theorem \ref{Theo4} bellow, we find many pairs $(N^{1},N^{2})$ of irreducible $\Free{+}$-modules such that $\Ext{\Free{+}}(N^2,N^1)=0$. 
In fact, we determine such pairs completely. 

By using the suitable pairs of irreducible $\Free{+}$-modules, we can prove that for irreducible $\charge{+}$-modules $M^{1}$ and $M^{2}$ in the fourth case, the exact sequence $0\rightarrow M^1\rightarrow M\rightarrow M^2\rightarrow0$ of $\charge{+}$-modules with a weak $V$-module $M$ splits as $\Free{+}$-modules.  
Thus $M=N^1\oplus N^2$ for some $\charge{+}$-submodule $N^1$ isomorphic to $M^1$ and $\Free{+}$-submodule $N^2$ isomorphic to $M^2$.
Then by using fusion rules for $\Free{+}$ determined in \cite{A1}, we can show that $N^2$ is in fact a $\charge{+}$-module and that $M$ is completely reducible as a $\charge{+}$-module.

This paper is organized as follows: 
In Section \ref{Prel} we recall definitions of vertex operator algebra, its module, Zhu's algebra, intertwining operator and fusion rule.
We also prepare categorical lemmas on short exact sequences for $V$-modules. 
The vertex operator algebras $\Free{+}$ and $\charge{+}$ and their irreducible modules are constructed in Section \ref{charge} following \cite{FLM}, \cite{D1} and \cite{D2}. 
We review the structure of Zhu's algebra $A(\Free{+})$ briefly following \cite{DN1} and describe the irreducible decompositions of irreducible $\charge{+}$-modules as $\Free{+}$-modules.
In Section \ref{vector}, we construct distinguished homogeneous vectors $H^{2r}$ in $\Free{+}$ for $r\in\Zpos$.
The grade-preserving action of these vectors are all diagonal on each irreducible $\Free{+}$-modules and the eigenvalues play a technical role to prove that the short exact sequence of $\Free{+}$-modules splits in Section \ref{split}. 
In Sections \ref{split} and \ref{indecomp}, we argue short exact sequences $0\rightarrow M^1\rightarrow M\rightarrow M^2\rightarrow0$ of $\Free{+}$-modules for irreducible $\Free{+}$-modules $M^1,M^2$ and a weak $\Free{+}$-module $M$.
All pairs $(M^1,M^2)$ of irreducible $\Free{+}$-modules for which the exact sequence splits are given in Section \ref{split}.
In Section \ref{indecomp}, for any pair $(M^1,M^2)$ of irreducible $\Free{+}$-modules not being indicated in Section \ref{split}, we construct indecomposable $\Free{+}$-modules $M$ for $\Free{+}$ such that the exact sequence $0\rightarrow M^1\rightarrow M\rightarrow M^2\rightarrow0$ does not split.
Section \ref{rationality} is devoted to a proof of rationality of the vertex operator algebra $\charge{+}$ for a rank one lattice $L$.
The result of Section \ref{rationality} extends to even lattices of arbitrary rank in Section \ref{Generalrank} with the help of Theorem 6.11 in \cite{Miya3}. 
In the Appendix \ref{appendix} we display the commutation relations related to the vectors $H^{2r}$ for $r=1,2,3$. 
  
{\bf Acknowledgments:} We would like to thank Masahiko Miyamoto and Chongying Dong for helpful comments and discussions.
We also thank Ching Hung Lam, Hiromichi Yamada, Atsushi Matsuo and Hiroki Shimakura for listening to my talk and giving me useful suggestions.

\section{Preliminaries}\label{Prel} 
In this section we review some definitions and notions from representation theory of vertex operator algebra (cf. \cite{B}, \cite{FLM}, \cite{FHL}, \cite{MN}, \cite{Zh1} and \cite{FZ}).

A vertex operator algebra $(V,Y(\,\cdot\,,z),\1,\w)$ is a $\Z$-graded vector space $V=\bigoplus_{n\in\Z}V_n$ with a linear map 
\[
Y(\,\cdot\,,z):V\to(\End V)[[z,z^{-1}]],\,a\mapsto Y(a,z)=Y(a,z)=\sum_{n\in\Z}a(n)z^{-n-1}
\] such that $\1\in V_0,\,\w\in V_2$ and that some axioms hold (see \cite{FLM}, \cite{FHL} and \cite{MN}). 
By definition, $V$ has distinguished vectors, the \textit{vacuum vector} $\1$ and the \textit{Virasoro element} $\w$ such that $Y(\1,z)=\id$ and $\{L(n):=\w(n+1),\id\}_{n\in\Z}$ gives a representation of the Virasoro algebra on $V$. 
Each $V_n\,(n\in\Z)$ is a finite dimensional eigenspace for $L(0)$ and $V_n=0$ for sufficiently small $n$.  

An \textit{automorphism} $g\in\End V$ of $V$ is a linear isomorphism satisfying $g(\w)=\w$ and $gY(a,z)g^{-1}=Y(g(a),z)$ for any $a\in V$. 
For any finite subgroup $G$ of automorphisms of $V$, the fixed point set $V^{G}=\{\,a\in V|g(a)=a\hbox{ for any }g\in G\,\}$ has naturally a vertex operator algebra structure. 
The Virasoro element and the vacuum vector of $V^G$ are given by $\w$ and $\1$ respectively. 

Let $g$ be an automorphism of $V$ of order $T$. 
Then $V$ decomposes into a direct sum of eigenspaces for $g$; $V=\bigoplus_{r=0}^{T-1}V^{r},\,V^r=\{\,a\in V|g(a)=e^{-2\pi ir/T}a\,\}$. 

A \textit{$g$-twisted weak $V$-module} $M$ is a vector space equipped with a linear map 
\[
Y(\,\cdot\,,z):V\to(\End M)[[z^{1/T},z^{-1/T}]]
\]
such that the following axioms hold on $M$ for any $a\in V^{r},\,b\in V$ and $u\in M$, 
\begin{align*}
Y(a,z)&=\sum_{n\in r/T+\Z}a(n)z^{-n-1}\quad (a(n)\in\End M),\\
a(n)u&=0\quad\hbox{for sufficiently large $n\in\Q$,}
\end{align*}

\noindent
\textit{(The twisted Jacobi identity)}
\begin{align}
\begin{split}\label{Jacobi}
&z_{0}^{-1}\delta\left({\frac{z_{1}-z_{2}}{z_{0}}}\right)Y_{M}(a,z_{1})Y_{M}(b,z_{2})-z_{0}^{-1}\delta\left({\frac{z_{2}-z_{1}}{ -z_{0}}}\right)Y_{M}(b,z_{2})Y_{M}(a,z_{1})\\
&\quad=z_{2}^{-1}\delta\left({\frac{z_{1}-z_{0}}{z_{2}}}\right)\left({\frac{z_{1}-z_{0}}{z_{2}}}\right)^{-{\frac{r}{T}}}Y_{M}(Y(a,z_{0})b,z_{2}),
\end{split}
\end{align}
and $Y_{M}(\1,z)=\id_{M}$. 
The identity \eqref{Jacobi} is equivalent to the Borcherds identity: For any $a,b\in V$ and $u\in M$ and $p\in\Z,s,t\in\Q$,
\begin{multline}\label{Borcherds}
\sum_{i=0}^{\infty}\binom{s}{i}(a(p+i)b)(s+t-i)u\\
=\sum_{i=0}^{\infty}(-1)^i\binom{p}{i}(a(p+s-i)b(t+i)-(-1)^pb(p+t-i)a(s+i))u.
\end{multline}
We also denote by $L(n)=\w(n+1)$ for any $n\in\Z$.
Then $\{L(n),\id\}$ gives $M$ a representation of the Virasoro algebra of central charge $c_V$.
Furthermore the $L(-1)$-derivative property 
\begin{align}\label{deriv1}
Y(L(-1)a,z)=\frac{d}{dz}Y(a,z)
\end{align}  
holds for any $a\in V$. 
A $g$-twisted weak $V$-module is called a \textit{weak $V$-module} in the case $g=\id$.
Such a module is often called an untwisted module.  

We shall recall some definitions related to modules only for untwisted modules for simplicity.   
A \textit{weak $V$-submodule} $N$ of a weak $V$-module $M$ is a subspace of $M$ such that $a(n)N\subset N$ for any $a\in V$ and $n\in\Z$. 
A \textit{homogeneous vector of weight} $r$ means an eigenvector for $L(0)$ of eigenvalue $r$.
We denote by $\wt{u}$ the weight of a homogeneous vector $u\in M$. 
We write $M_{(r)}$ for the generalized eigenspace for $L(0)$ of weight $r$ and $M_r$ for the eigenspace of weight $r$;
\begin{align*}
M_{(r)}&=\{\,u\in M\,|\,(L(0)-r)^ku=0\mbox{ for some}k\in\Z_{>0}\,\},\\
M_{r}&=\{\,u\in M\,|\,L(0)u=ru\,\}. 
\end{align*} 

A \textit{$\Zpos$-graded weak $V$-module} $M$ is a weak $V$-module which admits a $\Zpos$-grading $M=\bigoplus_{r\in\Zpos}M(r)$ such that the condition 
\begin{align}\label{admis}
a(n)M(m)\subset M(\wt{a}+m-n-1)
\end{align} 
is satisfied for any homogeneous $a\in V$ and $n\in\Z,m\in\Zpos$, where by definition $M(m)=0$ if $m<0$.
For a $\Zpos$-graded weak $V$-module $M$, we have the degree operator $d\in\End M$ defined by $d(u)=r u$ for any $u\in M(r)$. 
We call a subspace $N$ of $M$ a \textit{$\Zpos$-graded weak $V$-submodule} if $N$ is a $d$-stable weak $V$-submodule of $M$.   
A $\Zpos$-graded weak $V$-module $M$ is said to be \textit{irreducible} if there is no $\Zpos$-graded weak $V$-submodule except $0$ and $M$. 
We call a $\Zpos$-graded weak $V$-module $M$ \textit{completely reducible} if $M$ is a direct sum of irreducible $\Zpos$-graded weak $V$-submodules.
  
\begin{definition}
A vertex operator algebra $V$ is called \textit{rational} if any $\Zpos$-graded weak $V$-module is completely reducible.
\end{definition}

We call a weak $V$-module $M$ a \textit{$\C$-graded $V$-module} if $M$ has a $\C$-grading $M=\bigoplus_{r\in\C}M(r)$ satisfying the condition \eqref{admis} for any homogeneous $a\in V$ and $n\in\Z,m\in\C$ such that for any $r\in\C$, $M(r)$ is finite dimensional and $M_{r-n}=0$ for sufficiently large integer $n$.  
In the case $d=L(0)$, i.e., $M(r)=M_r$, $M$ is called a \textit{$V$-module}. 

Let $(M,Y(\,\cdot\,,z))$ be a $\C$-graded $V$-module. 
Then the restricted dual $M'=\bigoplus_{r\in\C}M(r)^{*}$ has a module structure $Y^*(\,\cdot\,,z)$ defined by the property
\[
\langle Y^*(a,z)v^{*},v\rangle=\langle v^{*},Y(e^{zL(1)}(-z^{-2})^{L(0)}a,z^{-1})v\rangle
\] 
for $a\in V,\,v\in M$ and $V^*\in M'$ (see \cite{FHL}), where $\langle\,\cdot,\cdot\,\rangle:M'\times M\to \C$ is the canonical pairing. 
The $\C$-graded $V$-module $(M',Y^*(\,\cdot\,,z))$ is called the \textit{contragredient module} to $M$. 
In general, the structure of contragredient module depends on the choice of gradings. 
When $M$ is a $\C$-graded $V$-module such that $M(r)=M_{(r)}$, i.e., $M$ is a direct sum of generalized eigenspaces for $L(0)$, then we always assume that the contragredient module is the restricted dual $M'=\bigoplus_{r\in\C}(M_{(r)})^{*}$ with the module structure $Y^*(\,\cdot\,,z)$.
As proved in \cite{FHL}, the contragredient module to an irreducible $\C$-graded $V$-module is also irreducible. 

We recall the notions of intertwining operator and fusion rule (see \cite{FHL}). 
For three weak $V$-modules $W^1,W^2$ and $W^3$, an \textit{intertwining operator of type $\fusion{W^1}{W^2}{W^3}$} is a linear map $\mathcal{Y}(\,\cdot\,,z):W^1\to\Hom(W^2,W^3)\{z\}$ which maps $v\in W^3$ to a formal power series $\mathcal{Y}(v,z)=\sum_{n\in\C}v(n)z^{-n-1}$ of $z$ with complex number power. 
By definition, for any $u\in W^1,\,v\in W^2$ and $h\in\C$, $u(h+n)v=0$ for sufficiently large integer $n$, and $a(n)$ and $b(m)$ for $a\in V,\,b\in W^1$ and $p,s\in\Z,\,t\in\C$ satisfies the Borcherds identity \eqref{Borcherds}. 
For any $u\in W^1$, the $L(-1)$-derivative property \eqref{deriv1} is also satisfied. 
The \textit{fusion rule} of type $\fusion{W^1}{W^2}{W^3}$ is the dimension of the vector space of all intertwining operator of type $\fusion{W^1}{W^2}{W^3}$. 

Let $M$ be a weak $V$-module. 
We introduce a shifted notation $\Tl{a}{n}$ for any $a\in V$ and $n\in\Z$ by setting 
\[
\Tl{a}{n}=a(\wt{a}+n-1)
\] 
for any homogeneous vector $a\in V$ and by extending it to $V$ linearly.
Then the Borcherds identity \eqref{Borcherds} leads the following two identities: The associativity formula 
\begin{align}\label{assohom}
\begin{split}
\Tl{(a(n)b)}{q}=&\sum_{i=0}^{\infty}\binom{n}{i}(-1)^{i}(\Tl{a}{-\wt{a}+1+n-i}\Tl{b}{\wt{a}-1-n+q+i}\\
&-(-1)^{n}\Tl{b}{\wt{a}-1+q-i}\Tl{a}{-\wt{a}+1+i}),
\end{split}
\end{align}
and the commutativity formula 
\begin{align}\label{commhom}
\begin{split}
[\Tl{a}{p},\Tl{b}{q}]=\sum_{i=0}^{\infty}\binom{\wt{a}+p-1}{i}\widetilde{(a(i)b)}(p+q)
\end{split}
\end{align}
for $a,b\in V$ and $n,p,q\in\Z$.
We note that if $M$ is a $\C$-graded module, then $\Tl{a}{n}M(m)\subset M(m-n)$ for any $a\in V,\,n\in\Z$ and $m\in\C$.

We next review the construction of Zhu's algebra introduced in \cite{Zh1} and recall its representation theory form \cite{Zh1} and \cite{DLM1}. 
Associated to a vertex operator algebra $V$, \textit{Zhu's algebra} $A(V)$ is defined as the quotient vector space $V/O(V)$, where $O(V)$ is the subspace spanned by vectors of the form $\sum_{i=0}^{\infty}\binom{\wt{a}}{i}a(i-2)b$ for any homogeneous vector $a\in V$ and $b\in V$. 
The product of $A(V)$ is induced from the binary operation $a*b=\sum_{j=0}^{\infty}\binom{\wt{a}}{j}a(j-1)b$ for any homogeneous vector $a\in V$ and $b\in V$. 
We denote by $[a]$ the image of $a\in V$ in $A(V)$.
Then $A(V)$ is an associative algebra with the product $[a]*[b]=[a*b]$ and has a unit $[\1]$ and a central element $[\w]$. 

For a weak $V$-module $M$, we set 
\[
\Omega(M)=\{\,u\in M\,|\,\Tl{a}{n}u=0\hbox{ for any $n\in \Z_{>0}$}\,\}.
\] 
Then $\Omega(M)$ is a representation space for $A(V)$ by the representation induced from the linear map $V\to\End\Omega(M),a\mapsto\Tl{a}{0}$. 
Conversely, in \cite{Zh1}, it is proved that for any $A(V)$-module $W$ there exists a $\Zpos$-graded weak $V$-module $M$ such that $\Omega(M)$ is isomorphic to $W$ as an $A(V)$-module. 
Furthermore, the map $M\mapsto \Omega(M)$ gives rise to a one to one correspondence between the set of equivalence classes of irreducible $\Zpos$-graded weak $V$-modules and that of equivalence classes of irreducible $A(V)$-modules. 
It is proved in \cite{Zh1} and \cite{DLM1} that $V$ is rational then Zhu's algebra $A(V)$ is semisimple, that there are only finitely many irreducible $\Zpos$-graded weak $V$-modules up to isomorphisms and that any irreducible $\Zpos$-graded weak $V$-module is a $V$-module. 

Let $C_1(V)$ be the subspace of $V$ spanned by the subspace $L(-1)V$ and the set of vectors $a(-1)b$ for $a,b\in\bigoplus_{d=1}^{\infty}V(d)$. 
Let $S$ be a set of homogeneous vectors in $V$ such that $V=\textrm{span}\,S+C_1(V)$.
Then it is proved in \cite{KL} that $V$ is spanned by vectors of the form $a_{1}(-n_{1})\ldots a_{s}(-n_{s})\1$ for $s\geq0,\,a_{i}\in S$ and $n_{i}\geq1$. 
So we call such a subset $S$ a \textit{set of generators of $V$} .
\begin{proposition}\label{aaee}
Suppose that $V=\bigoplus_{n=0}^{\infty}V_n$ and that $V_0=\C\1$. 
Let $M$ be a weak $V$-module and $S$ a set of generators of $V$. 
Then $u\in\Omega(M)$ if and only if $\Tl{a}{n}u=0$ for any $a\in S$ and $n\geq1$. 
\end{proposition}
\begin{proof}
It is clear that if $u\in\Omega(M)$ then $\Tl{a}{n}u=0$ for any $a\in S$ and $n\geq1$.
Conversely, let $u\in M$ and suppose $\Tl{a}{n}u=0$ for any $a\in S$ and $n\geq1$.  
Set $X=\{a\in V|\,\Tl{a}{n}u=0\hbox{ for any $n\geq1$}\}$.   
By using induction on the weight of a homogeneous vector $a\in V$, we shall prove that $a\in X$.  
It is trivial that $V_0=\C\1\subset X$. 
Assume that $\bigoplus_{n=0}^{r}V_n\subset X$ for some $r\geq0$. 
Then we have to prove that $a=a_{1}(-n_{1})\ldots a_{s}(-n_{s})\1\in X$ with $\wt{a}=r+1$ for $a_{i}\in S$ and $n_{i}\geq1$.
We may assume that $\wt{a_1}>0$. 
Then $v:=a_{2}(-n_{2})\ldots a_{s}(-n_{s})\1\in X$ because $\wt{v}\leq r$.   
We also see that $\wt{a_1(j)v}<\wt{a}$ for any $j\in\Zpos$ and hence $a_1(j)v\in X$ by induction hypothesis. 
Thus the commutativity formula \eqref{commhom} shows that $\Tl{a_1}{p}\Tl{v}{q}u=\Tl{v}{q}\Tl{a_1}{p}u=0$ for any $p,q\in\Z$ with $p+q\geq1$. 
Therefore, the associativity formula \eqref{assohom} implies that 
\begin{align*}
\Tl{(a_1(-n_1)v)}{q}u=&\sum_{i=0}^{\infty}\binom{-n_1}{i}(-1)^{i}(\Tl{a_1}{-\wt{a_1}+1-n_1-i}\Tl{v}{\wt{a_1}-1+n_1+q+i}u\\
&-(-1)^{n_1}\Tl{v}{\wt{a_1}-1+q-i}\Tl{a_1}{-\wt{a_1}+1+i})u=0
\end{align*}
if $q\geq1$.  
This proves that $a_1(-n_1)v\in X$.
\end{proof}

We finally state some lemmas on short exact sequences of $V$-modules.  
We first give the following lemma whose proof is elemental. 
\begin{lemma}\label{Prop53}
Let $0\rightarrow M^{1}\rightarrow M\rightarrow M^2\rightarrow 0$ be a short exact sequence of weak $V$-modules and  $a\in V$.
Suppose that $\Tl{a}{0}$ acts on both $M^1$ and $M^2$ diagonally. 
Then $M$ is a direct sum of generalized eigenspaces for $\Tl{a}{0}$. 
Furthermore, the eigenspace of $M^2$ of eigenvalue $h$ is the image of the generalized eigenspace of $M$ of eigenvalue $h$. 
\end{lemma}

We note that if $a$ and $b\in V$ satisfy the condition in Lemma \ref{Prop53} and if $\Tl{a}{0}$ and $\Tl{b}{0}$ commute, then $M$ is a direct sum of simultaneous generalized eigenspaces for $\Tl{a}{0}$ and $\Tl{b}{0}$. 

Let $M^1$ and $M^2$ be weak $V$-modules.
Following \cite{MP} we define the extension group $\Ext{V}(M^2,M^1)$ of $M^2$ by $M^1$. 
We call a weak $V$-module $M$ an \textit{extension of $M^2$ by $M^1$} if there is a short exact sequence $0\rightarrow M^{1}\rightarrow M\rightarrow M^2\rightarrow 0$. 
Two extensions $M$ and $N$ of $M^2$ by $M^1$ are said to be equivalent if there exists a $V$-module homomorphism $f:M\to N$ such that the following diagram commutes:
\begin{align*}
\begin{CD}
0@>>> M^1 @>{\phi}>>  M @>{\psi}>> M^2@>>>0\quad\mbox{(exact)}\\
 @.       @|            @VV{f}V       @|  \\
0@>>> M^1 @>{\phi'}>> N @>{\psi'}>>M^2@>>>0\quad\mbox{(exact)}\\
\end{CD}
\end{align*}
We define $\Ext{V}(M^2,M^1)$ the set of all equivalent classes of extensions of $M^2$ by $M^1$.
The set $\Ext{V}(M^2,M^1)$ has a structure of an abelian group such that the equivalence class of $M^1\oplus M^2$ is zero (see \cite[Section 1.12]{MP}).  

We next give a sufficient condition to show that $\Ext{V}(M^2,M^1)=0$. 
\begin{proposition}\label{Prop73}
Let $M^1$ and $M^2$ be $V$-modules.
Suppose that the sets of weights on $M^1$ and $M^2$ have no intersection.
Then $\Ext{V}(M^2,M^1)=0$. 
\end{proposition}
\begin{proof}
Let $0\rightarrow M^{1}\rightarrow M\xrightarrow{\psi} M^2\rightarrow 0$ be a short exact sequence. 
By using Lemma \ref{Prop53} with $a=\w$, we see that $M$ is a $\C$-graded $V$-module by the generalized eigenspace decomposition for $L(0)$.
Let $P_i$ be the set of weights on $M^i$ for $i=1,2$. 
Then $M=(\bigoplus_{r\in P^1}M_{(r)})\oplus(\bigoplus_{r\in P^2}M_{(r)})$ by the assumption.
Then it is clear that $M^1\cong\Ker \psi=\bigoplus_{r\in P^1}M_{(r)}$ and that the restriction of $\psi$ to $\bigoplus_{r\in P^2}M_{(r)}$ induces an isomorphism from $\bigoplus_{r\in P^2}M_{(r)}$ to $M^2$. 
This proves $M\cong M^1\oplus M^2$. 
\end{proof}
 
\begin{proposition}\label{Prop7}
Let $M$ and $N$ be irreducible $V$-modules. 
Then $\Ext{V}(N,M)=0$ if and only if $\Ext{V}(M',N')=0$.
\end{proposition}
\begin{proof}
Since $M$ and $N$ are $V$-modules, $(M')'\cong M$ and $(N')'\cong N$ as $V$-modules.
Thus it suffices to show that if $\Ext{V}(M',N')=0$ then $\Ext{V}(N,M)=0$.

Let $0\rightarrow M\,\xrightarrow{\phi}\,L\,\xrightarrow{\xi}\, N\rightarrow 0$ be a short exact sequence or a weak $V$-module $L$.  
By Lemma \ref{Prop53}, $L=\bigoplus_{r\in\C}L_{(r)}$ and $\dim L_{(r)}$ is finite for any $r\in\C$.
Thus $L$ is a $\C$-graded $V$-module. 
Then by taking the contragredient module to $L$ with respect to this grading, we have the exact sequence $0\rightarrow N'\,\xrightarrow{\xi^*}\, L'\,\xrightarrow{\phi^*}\, M'$, where the $V$-module homomorphism $\phi^*$ is given by $\phi^*(f)=f\circ\phi$ for any $f\in L'$ and that the map $\xi^*$ is given by the same way. 
Since $\phi^*$ is nonzero and $M'$ is irreducible, $\phi^*$ is surjective.
Thus $L'\cong M'\oplus N'$ as $V$-modules because $\Ext{V}(N',M')=0$.
Noting that $(L_{(r)})^*=(L')_{(r)}$ for any $r\in\C$ and that $(L')'\cong L$, we see that $L\cong N\oplus M$.
\end{proof}
\begin{lemma}\label{lemma3}
Let $\{M^{i}\}_{i\in I}$ and $\{N^{j}\}_{j\in J}$ be sets of irreducible $V$-modules.
Suppose that $\Ext{V}(N^{j},M^{i})=0$ for any $i\in I,\,j\in J$ and that $\bigoplus_{i\in I}M^{i}$ is a $V$-module. 
Then $\Ext{V}(\bigoplus_{j\in J}N^{j},\bigoplus_{i\in I}M^{i})=0$. 
\end{lemma}
\begin{proof}
We first consider the case $J=\{j_0\}$ is a one point set, and set $N=N^{j_0}$. 
Let $0\rightarrow \bigoplus_{i\in I}M^{i}\,\xrightarrow{\phi}\, M\,\xrightarrow{\xi}\, N\rightarrow 0$ be a short exact sequence for a weak $V$-module $M$.  
We denote by $\phi_{i}$ the restriction to $M^{i}$ of the $V$-module homomorphism $\phi$.
For any $i\in I$, we set $K^i=\bigoplus_{j\in I-\{i\}}\phi(M^{j})\subset M$ and consider the quotient module $M/K^i$.
Let $\pi_{i}$ and $\bar{\xi}$ be the canonical projection $M\to M/K^i$ and the induced surjection $M/K^i\to N$ respectively. 
Then we have a short exact sequence 
\[
0\rightarrow M^{i}\xrightarrow{\pi_i\circ \phi_i} M/K^i\xrightarrow{\bar{\xi}} N\rightarrow 0.
\]
Therefore, by the assumption of the lemma, there exists a $V$-module homomorphism $\psi_{i}$ from $M/K^i$ to $M^{i}$ such that $\psi_{i}\circ(\pi_{i}\circ \phi_{i})=\id_{M^{i}}$.  
Therefore, the composition $\psi_{i}\circ\pi_{i}:M\to M^{i}$ satisfies the property $(\psi_{i}\circ\pi_{i})\circ \phi_{i}=\id_{M^{i}}$.
We define a $V$-module homomorphism $\psi:M\to \bigoplus_{j\in I}M^{j}$ by $\psi(u)=(\psi_{i}\circ\pi_{i}(u))_{i\in I}$. 
If $\psi$ is well defined, then $\psi\circ\phi=\id$, and this proves the lemma.

To prove that $\psi$ is well defined, we need to show that for any $u\in M$, $\psi_{i}\circ\pi_{i}(u)=0$ except for finite many $i$. 
Since $M=\bigoplus_{r\in\C}M_{(r)}$ by Lemma \ref{Prop53}, it suffices to show the claim that for any $r\in\C$, $\psi_i\circ\pi_i(M_{(r)})=0$ except for finitely many $i$.
It is easy to see that $\psi_i\circ\pi_i(M_{(r)})\subset (M^i)_r$. 
On the other hand, there are only finite many indices $j$ so that $(M^j)_r\neq0$ for given $r\in \C$ by the assumption. 
This proves the claim. 

In the case $J$ is a general set, let
\begin{align}\label{exa31}
0\rightarrow \bigoplus _{i\in I}M^{i}\xrightarrow{\phi} M\xrightarrow{\xi} \bigoplus_{j\in J}N^j\rightarrow 0
\end{align}
be a short exact sequence. 
Then we have the short exact sequence $0\rightarrow \bigoplus _{i\in I}M^{i}\xrightarrow{\phi} \xi^{-1}(N^j)\xrightarrow{\xi}N^j\rightarrow 0$ for any $j\in J$ and this exact sequence splits.
Hence there exists $s_{j}:N^j\to \xi^{-1}(N^j)$ such that $\xi\circ s_j=\id$.
Then $s=\oplus_{j\in J}s_j$ satisfies $\xi\circ s=\id$.  
This shows that the exact sequence \eqref{exa31} splits.
\end{proof}

We recall the notion of $C_2$-cofiniteness. 
Set $C_2(V):=\haru{a(-2)b}{a,b\in V}$. 
Then a vertex operator algebra $V$ is called \textit{$C_2$-cofinite} if $C_2(V)$ is of finite codimension in $V$. 
A $C_2$-cofinite vertex operator algebra has the following remarkable properties (see \cite{Miya}).
\begin{theorem}\label{C2-condition} 
Let $V$ be a $C_2$-cofinite vertex operator algebra.
Then any weak $V$-module $M$ is decomposed into a direct sum of generalized eigenspaces for $L(0)$, that is, $M=\bigoplus_{r\in \C}M_{(r)}$.
For any nonzero $w\in M_{(r)}\,(r\in\C)$, the weak $V$-submodule generated from $w$ is $\C$-graded.  
\end{theorem}
This theorem implies that if $V$ is $C_2$-cofinite then weak $V$-module is a sum of $\C$-graded $V$-submodules. 
In particular we have the following theorem:
\begin{theorem}\label{extension}
Let $V$ be a $C_2$-cofinite vertex operator algebra.
Then $V$ is rational if and only if $\Ext{V}(N,M)=0$ for any pair of irreducible $V$-modules $M$ and $N$.  
\end{theorem}
\begin{proof}
Suppose that $V$ is rational.
Let $M$ and $N$ be irreducible $V$-modules and consider a short exact sequence $0\rightarrow M\rightarrow L\rightarrow N\rightarrow0$ for a weak $V$-module. 
Because $L$ is $\C$-graded, $L$ is completely reducible. 
Thus $L\cong M\oplus N$.  

Conversely, we suppose that $\Ext{V}(N,M)=0$ for any pair of irreducible $V$-modules $M$ and $N$.
By Theorem \ref{C2-condition}, it suffices to show that any $\C$-graded $V$-module $M=\bigoplus_{r\in \C}M_{(r)}$ is completely reducible. 
Since $V$ is $C_2$-finite, $A(V)$ is finite dimensional; we note that the identity of $V$ induces a surjective map from $V/C_2(V)$ to ${\rm gr}A(V)$, where ${\rm gr}A(V)$ is the graded algebra of the filtered algebra $A(V)$ with filtration by weights (see \cite{Zh1}).
Thus the $A(V)$-module $\Omega(M)$ has an irreducible submodule. 
We take an irreducible $A(V)$-submodule $W$ such that the weight $r$ for $L(0)$ is maximal in the sense that $s-r$ is not a positive integer for any eigenvalue for $L(0)$ on $\Omega(M)$.   
Then $W$ generates an irreducible $V$-submodule of $M$. 
This shows that $M$ has an irreducible $V$-submodule. 

Let $N$ be the sum of all irreducible $V$-submodules of $M$.
Then it is a $\Zpos$-graded weak $V$-submodule of $M$. 
Assume that $M\neq N$. 
Consider the quotient module $M/N$ and its irreducible submodule $M^0/N$, where $M^0$ is a $\Zpos$-graded weak $V$-submodule containing $N$. 
Then there is a generalized eigenvector $w\in M^0$ for $L(0)$ such that $M^{0}=\langle\,w\,\rangle+N$, where $\langle\,w\,\rangle$ is the $\Zpos$-graded weak $V$-submodule of $M^0$ generated from $w$. 
By Theorem \ref{C2-condition}, $\langle\,w\,\rangle$ is a $\C$-graded $V$-module. 
Therefore, $\langle\,w\,\rangle\cap N$ is a direct sum of irreducible $V$-modules and $\langle\,w\,\rangle/(\langle\,w\,\rangle\cap N)\cong M^0/N$ is irreducible. 
By Lemma \ref{lemma3}, we have $\Ext{V}(M^0/N,\langle\,w\,\rangle\cap N)=0$ and hence $M^0\cong N\oplus M^0/N$ as a $V$-module. 
This is a contradiction. 
Hence $M$ is completely reducible. 
\end{proof}
\section{Vertex operator algebras $\Free{+}$ and $\charge{+}$}\label{charge}
In this section we shall review a construction of vertex operator algebras $\Free{+}$ and $\charge{+}$ and recall the construction and the classification of their irreducible modules. 
We also state irreducible decomposition of irreducible $\charge{+}$-modules as $\Free{+}$-modules and the structure of Zhu's algebra $A(\Free{+})$ in the rank one case.  

Let $L$ be an even lattice with positive definite $\Z$-bilinear form $(\,\cdot\,,\cdot\,)$. 
Set $\h=\C\otimes_{\Z} L$ and extend $(\,\cdot\,,\cdot\,)$ to the nondegenerate symmetric $\C$-bilinear form $(\,\cdot\,,\cdot\,)$.
Then $\hh=\h\otimes\C[x,x^{-1}]\bigoplus\C C$ is a Lie algebra by the commutation relation
\begin{align}\label{commutation}
[X_{1}\otimes x^{m},X_{2}\otimes x^{n}]=m(X_1,X_2)\delta_{m,-n}C,\quad [C,\hh]=0
\end{align}
for $X_{1},X_{2}\in\h$ and $m,n\in\Z$. 
We see that $\hh^{+}=\h\otimes\C[x]\bigoplus\C C$ is a commutative Lie subalgebra of $\hh$.
For $\lambda\in\h$, we define an $\hh^{+}$-module structure on the one dimensional vector space $\C e^{\lambda}$ by the action 
\[
\rho(X\otimes x^{n})e^{\lambda}=(h,\lambda)\delta_{n,0}e^{\lambda}\quad\hbox{for $X\in\h,\,n\geq0$.}
\]
We denote by $\Fremo{\lambda}$ the induced module $U(\hh)\otimes_{U(\hh^{+})}\C e^{\lambda}$, where $U(\g)$ is the universal enveloping algebra of $\g$.  
The action of $X\otimes t^{n}$ for $X\in\h,\,n\in\Z$ will be written by $X(n)$. 

We set $\Free{}=\Fremo{0}$ and $\1=1\otimes e^{0},\,\w=\frac{1}{2}\sum_{i=0}^{d}h(-1)^2\1$, where $d=\dim\h$ and $\{h_i\}$ is an orthonormal basis of $\h$. 
Then there exists a linear map $Y(\,\cdot\,,z):\Free{}\to(\End\Fremo{\lambda})[[z,z^{-1}]]$ for $\lambda\in\h$ such that $(\Free{},Y(\,\cdot\,,z),\1,\w)$ becomes a simple vertex operator algebra and $(\Fremo{\lambda},Y(\,\cdot\,,z))$ is an irreducible $\Free{+}$-module (see \cite{FLM}). 

For a vector $a=X_{1}(-n_{1})\cdots X_{s}(-n_{s})\1\,(X_i\in\h,\,n_j\geq1)$, the action of $\Tl{a}{n}$ for $n\in\Z$ on $\Fremo{\lambda}$ are given by the formula 
\begin{align}\label{action34}
\Tl{a}{n}v=\sum_{\substack{i_j\in\Z,\, \sum i_j=n}}\binom{-i_1-1}{n_1-1}\cdots\binom{-i_s-1}{n_s-1}\NO X_1(i_1)\cdots X_s(i_s)\NO v
\end{align}
for $v\in\Fremo{\lambda}$, where the normal ordering product $\NO\cdot\NO$ is the operation to reorder $X(n)$ with $n\geq1$ to the right of $X'(n)$ with $n\leq0$. 
We see that the weight of the vector $X_{1}(-n_{1})\cdots X_{s}(-n_{s})e^{\lambda}\in\Fremo{\lambda}$ is $\sum_{i=1}^{s}n_i+\frac{(\lambda,\lambda)}{2}$. 

Set $\C[\h]=\bigoplus_{\lambda\in \h}\C e^{\lambda}$.
For a subset $M\in\h$, let $\C[M]=\bigoplus_{\lambda\in M}\C e^{\lambda}$. 
We note that $\C[M]$ has an $\hh^+$-module structure. 
Define $V_M=U(\hh)\otimes_{U(\hh^+)}\C[M]$.
By definition, $V_M$ has naturally an $\Free{}$-module structure such that $V_M\cong\bigoplus_{\lambda\in M}\Fremo{\lambda}$ as an $\Free{}$-module. 
We see that $\Free{}\subset \charge{}$ and hence $\1,\w\in \charge{}$.
Let $L^{\circ}=\{\lambda\in\h\,|\,(\lambda,L)\subset\Z\,\}$ be the dual lattice of $L$.  
In \cite{FLM} it is proved that for any $\lambda\in L^\circ$, there exists a linear map $Y(\,\cdot\,,z):\charge{}\to\End \charlam{\lambda}[[z,z^{-1}]]$ extending the module structure for $\Free{}$ such that $(\charge{},Y(\,\cdot\,,z),\1,\w)$ is a vertex operator algebra for $\lambda\in L$ and that $(\charlam{\lambda},Y(\,\cdot\,,z))$ is a $\charge{}$-module for any $\lambda\in L^{\circ}$.

We define a linear map $\theta:V_{\h}\to V_{\h}$ by 
\begin{align}\label{acttheta1}
\theta(X_{1}(n_{1})\cdots X_{r}(n_{r})e^{\lambda})=(-1)^{r} X_{1}(n_{1})\cdots X_{r}(n_{r})e^{-\lambda}
\end{align}
for $\lambda,h_i\in\h$ and $n_{i}\in\Z$.
By abuse of notation we also use $\theta$ when $\theta$ is restricted to a subset of $V_{\h}$.  
We write $W^\pm=\{u\in W\,|\,\theta(u)=\pm u\,\}$ respectively for a $\theta$-stable subset $W\subset V_{\h}$. 
The linear map $\theta$ gives an automorphism of $\charge{}$ and $\Free{}$ of order $2$, and hence $\charge{+}$ and $\Free{+}$ has naturally vertex operator algebra structures.
The subset $\Free{-}$ becomes an irreducible $\Free{+}$-modules and $\Fremo{\lambda}$ is an irreducible $\Free{+}$-module for $\lambda\in\h-\{0\}$.
We see that $\theta:\Fremo{\lambda}\rightarrow\Fremo{-\lambda}$ is an $\Free{+}$-module isomorphism for $\lambda\in\h$.
As for $\charge{+}$, $\charge{-}$ is an irreducible $\charge{+}$-module and $\charlam{\lambda}$ is an irreducible $\charge{+}$-module if $\lambda\in L^\circ-\frac{1}{2}L$. 
In the case $\lambda\in L^\circ\cap\frac{1}{2}L$, then $\charlam{\lambda}^{\pm}$ are irreducible $\charge{+}$-modules (see \cite{DLi}, \cite{DM}).
The linear map $\theta:\charlam{\lambda}\rightarrow\charlam{-\lambda}$ is a $\charge{+}$-module isomorphism for any $\lambda\in L^\circ$.   

Next we review a construction of the $\theta$-twisted modules for $\Free{}$ and $\charge{}$. 
Let $\hh[-1]=\h\otimes x^{\frac{1}{2}}\C[x,x^{-1}]\oplus\C C$ be a Lie algebra with commutation relation \eqref{commutation} for $X_1,X_2\in\h$ and $m,n\in\frac{1}{2}+\Z$. 
Set $\hh[-1]^{+}=\h\otimes x^{\frac{1}{2}}\C[x]\oplus\C C$, and consider the one dimensional $\hh[-1]^+$-module $\C\vactw$ with actions $\rho(X\otimes x^{n})\vactw=0,\,C\vactw=\vactw$ for $X\in\h$ and $n\in\frac{1}{2}+\Zpos$. 
We write $\Fretw{}$ for the induced module $U(\hh[-1])\otimes_{U(\hh[-1]^{+})}\C\vactw$.  
The action of $h\otimes x^{n}$ for $h\in\h,\,n\in\frac{1}{2}+\Z$ is also denoted by $h(n)$.
In \cite{FLM}, it is prove that there is a linear map $Y(\,\cdot\,,z):\Free{}\to(\End\Fretw{})[[z,z^{-1}]]$ such that $(\Fretw{},Y(\,\cdot\,,z))$ becomes a $\theta$-twisted $\Free{}$-module. 
It is known that $\Fretw{}$ is the unique  $\theta$-twisted $\Free{}$-module up to isomorphisms. 

Any irreducible $\theta$-twisted $\charge{}$-module is realized as a tensor product of $\Fretw{}$ and an irreducible module for a suitable group whose action commutes with that of $\Free{}$ (see \cite{FLM} and \cite{D2}). 
In the case of rank one, let $T_i$ be an irreducible $L/2L$-module on which $\alpha+2L$ acts by the scalar $(-1)^i$ for $i=1,2$.
Then $\charge{T_i}=\Fretw{}\otimes T_i$ has a $\theta$-twisted $\charge{}$-module structure.
In \cite{FLM}, it is proved that there exists a linear map $Y(\,\cdot\,,z):\charge{}\to\End (\charge{T_i})[[z^{1/2},z^{-1/2}]]$ extending the module structure for $\Free{}$ such that $(\charge{T_i},Y(\,\cdot\,,z))$ becomes a $\theta$-twisted $\charge{}$-module. 
Since $T_i$ is one dimensional, we have $\charge{T_i}\cong\Fretw{}$ as a $\theta$-twisted $\Free{}$-module.  

Now we define the action of $\theta$ on $\Fretw{}$ by 
\begin{align}\label{acttheta2}
\theta(X_{1}(-n_{1})\cdots X_{r}(-n_{r})\vactw)=(-1)^{r}X_1(-n_{1})\cdots X_r(-n_{r})\vactw
\end{align} 
for $h_i\in\h$ and $n_{i}\in\frac{1}{2}+\Zpos$.
Then the $\pm1$-eigenspaces $\Fretw{\pm}$ for $\theta$ are irreducible $\Free{+}$-modules. 
The action of $\theta$ is extended to $\charge{T_i}$ so that $\theta(u\otimes t)=\theta(u)\otimes t$ for $u\in\Fretw{}$ and $t\in T_i$. 
Then $\charge{T_i,\pm}=\{\,u\in\charge{T_i}\,|\,\theta(u)=\pm u\,\}$ are irreducible $\charge{+}$-modules (see \cite{DLi} and \cite{DN1}).

The irreducible modules for $\Free{+}$ and $\charge{+}$ are classified in \cite{DN1}--\cite{DN3} and \cite{AD}. 
We here state the result only for the rank one case. 
\begin{theorem} (1) {\rm (\cite{DN1})}
The set 
\[
\{\Free{\pm},\,\Fretw{\pm},\,\Fremo{\lambda}(\cong\Fremo{-\lambda})\,|\,\lambda\in\h-\{0\}\,\}
\]
gives the set of all inequivalent irreducible $\Free{+}$-modules.
 
(2) {\rm (\cite{DN2})}
Let $L=\Z\alpha$. 
The set 
\[
\{\charge{\pm},\,\charlam{\frac{\alpha}{2}}^{\pm},\,\charlam{\frac{r\alpha}{(\alpha,\alpha)}},\,\charge{T_1,\pm},\,\charge{T_2,\pm}\,|\,1\leq r\leq (\alpha,\alpha)/2-1\,\}
\]
gives the set of all inequivalent irreducible $\charge{+}$-modules.
\end{theorem}

Let $L=\Z\alpha$ be a positive definite even lattice of rank one.  
Then all irreducible $\charge{+}$-modules are decomposed into direct sums of irreducible $\Free{+}$-modules as follows (see \cite{DG} and \cite[Proposition 3.2]{A3}):
\begin{align}
\charge{\pm}&\cong\Free{\pm}\oplus\bigoplus_{m=1}^{\infty}\Fremo{m\alpha},\label{dec1}\\
\charlam{\lambda}&\cong\bigoplus_{m\in\Z}\Fremo{\lambda+m\alpha}\quad\hbox{for $\lambda\in L^{\circ}$},\label{dec2}\\ 
\charlam{\frac{\alpha}{2}}^{\pm}&\cong\bigoplus_{m=0}^{\infty}\Fremo{{\alpha}/{2}+m\alpha},\label{dec3}\\
\charge{T_{i},\pm}&\cong\Fretw{\pm}\quad\hbox{for $i=1,2$.}\label{dec4}
\end{align}
We see that a multiplicity of irreducible $\Free{+}$-module in an irreducible $\charge{+}$-module is at most one. 
This fact will be used in the later section. 

We finally recall the structure of Zhu's algebra $A(\Free{+})$ following \cite{DN1}.
Set $J=h(-1)^4\1-3h(-3)h(-1)\1+\frac{3}{2}h(-2)^2\1\in\Free{+}$. 
The vector $J$ is a singular vector for the Virasoro algebra of weight $4$. 
Then the set $\{\w,J\}$ is a set of generators of $\Free{+}$ (see \cite[Lemma 3.1]{DN1}). 
Furthermore, the images $[\w]$ and $[J]$ in $A(\Free{+})$ generate Zhu's algebra $A(\Free{+})$ with the following defining relations (see \cite[Theorem 4.4]{DN1}) 
\begin{align}
&[\w]*[J]-[J]*[\w]=0,\label{relations1}\\
&([\w]-1)*([\w]-1/16)*([\w]-9/16)([J]-[\w]+4[\w]*[\w])=0,\label{relations2}\\
&([J]-[\w]+4[\w]*[\w])(70 [J]+908[\w]*[\w]-515[\w]+27)=0.\label{relations3}
\end{align}

\section{The vectors $H^{2r}$ with $r\geq1$}\label{vector}
In Sections \ref{vector}--\ref{indecomp}, we develop the representation theory for $\Free{+}$ from the view of extensions of irreducible modules in the case of rank one. 
We will find pairs $(M^1,M^2)$ so that $\Ext{\Free{+}}(M^2,M^1)=0$ for irreducible $\Free{+}$-modules in Section \ref{split} and prove that the other pairs of irreducible $\Free{+}$-modules give nontrivial extensions in Section \ref{indecomp}. 
To demonstrate this, we will use distinguished vectors $H^{4}$ and $H^6$ in $\Free{+}$ (see \eqref{kkkk1} and \eqref{kkkk2}), whose grade-preserving actions $\Tl{H^4}{0}$ and $\Tl{H^6}{0}$ are diagonal on all irreducible $\Free{+}$-modules.
In this section, we construct more generally the series of vectors $H^{2r}$ for $r\geq1$ such that $H^{2}=\w$ and give some commutation relation among the modes of $\w,H^4,H^6$.  
  
We set $\mathcal{S}=\haru{\1,h(-n)h(-m)\1}{n,m\geq1}\subset\Free{+}$.
Then we can prove that $a(i)b\in \mathcal{S}$ for any $a,b\in \mathcal{S}$ and $i\in\Zpos$. 
In particular, $\mathcal{S}$ contains the Virasoro element $\w$ and is closed under the action of $L(-1)=\w(0)$. 
We see that $\mathcal{S}=\bigoplus_{r=0}^{\infty}\mathcal{S}_r,\mathcal{S}_r=\mathcal{S}\cap V_r$, and have $\mathcal{S}_0=\C\1,\,\mathcal{S}_1=0$ and $\mathcal{S}_2=\C\w$.

Now we consider the vector space $\mathcal{S}_r$ for any $r\geq2$.
By using the identity 
\[
L(-1)h(-n)h(-m)\1=nh(-n-1)h(-m)\1+mh(-n)h(-m-1)\1
\] 
for any $n,m\in\Zpos$, one can easily see that 
\begin{align}\label{SR1}
\mathcal{S}_r=
\begin{cases}L(-1)\mathcal{S}_{r-1}&\hbox{if $r$ is odd},\\
\C h(-r/2)^2\1\oplus L(-1)^2\mathcal{S}_{r-2}&\hbox{if $r$ is even.}
\end{cases}
\end{align}
\begin{lemma}\label{Lemma17}
For any $r\in\Z_{>0}$, there exists $H^{2r}\in \mathcal{S}_{2r}$ uniquely such that 
\[
[\Tl{H^{2r}}{0},h(n)]=-n^{2r-1}h(n)
\] 
for $n\in\frac{1}{2}\Z$ on any weak ($\theta$-twisted) $\Free{}$-module. 
Moreover, $\mathcal{S}_{2r}=\C H^{2r}\oplus L(-1)^2\mathcal{S}_{2r-2}$. 
\end{lemma}
\begin{proof}
We use induction on $r$. 
The case $r=1$, $H^{2}=\w$ satisfies the desired property and it is unique.  
Suppose that there exist $H^{2i}$ uniquely for $1\leq i\leq r-1$ satisfying the conditions in the lemma.
We calculate the commutation relation $[\Tl{h(-r)^2\1}{0},h(n)]$. 
Since $h(i)h(-r)^2\1=2r\delta_{i,r}h(-r)\1=\frac{2r}{(r-1)!}\delta_{i,r}L(-1)^{r-1}h(-1)\1$ by \eqref{action34}, we see that 
\begin{align}\label{eqn-2}
\begin{split}
[h(n), \Tl{(h(-r)^2\1)}{0}]&=\frac{2r}{(r-1)!}\binom{n}{r}\Tl{(L(-1)^{r-1}h(-1)\1)}{n}\\
&=(-1)^{r-1}2r\binom{n+r-1}{r-1}\binom{n}{r}h(n)
\end{split}
\end{align}
for any $n\in\frac{1}{2}\Z$.
We can find that $(-1)^{r-1}2r\binom{n+r-1}{r-1}\binom{n}{r}=\sum_{i=1}^{r}c_{i}n^{2i-1}$ for some $c_i\in\Q$ with $c_r\neq0$. 
Define a vector $H^{2r}$ by using $H^{i}$ for $1\leq i\leq r-1$ as follows:
\[
H^{2r}:=\frac{1}{c_r}\left(h(-r)^2\1-\sum_{i=1}^{r-1}\frac{c_i(2i-1)!}{(2r-1)!}L(-1)^{2r-2i}H^{2i}\right).
\]
Then $\Tl{H^{2r}}{0}=\frac{1}{c_r}\left(\Tl{(h(-r)^2\1)}{0}-\sum_{i=1}^{r-1}c_i\Tl{H^{2i}}{0}\right)$.
Hence by \eqref{eqn-2} and induction hypothesis, we get
\begin{align*}
[h(n), c_r \Tl{H^{2r}}{0}]&=[h(n),\Tl{(h(-r)^2\1)}{0}]-\sum_{i=1}^{r-1}c_i[h(n),\Tl{H^{2i}}{0}]\\
&=\sum_{i=1}^{r}c_in^{2i-1}h(n)-\sum_{i=1}^{r-1}c_i n^{2i-1}h(n)\\
&=c_r n^{2r-1}h(n)
\end{align*}
for any $n\in\frac{1}{2}\Z$. 
Therefore, the commutation relation $[\Tl{H^{2r}}{0},h(n)]=-n^{2r-1}h(n)$ holds for any $n\in\frac{1}{2}\Z$. 
From the construction of $H^{2r}$, the uniqueness of $H^{2r}$ and the second assertion are clear.
\end{proof}
\begin{remark}\label{remark2}
Let $c_i\,(1\leq i\leq r)$ be the constant in the proof of Lemma \ref{Lemma17}.
By direct calculations, we see that $c_1=2$ and $\Tl{(h(-r)^2\1)}{0}e^{\lambda}=(\lambda,\lambda)e^{\lambda}=2L(0)e^{\lambda}$ for any $\lambda\in\h$.
Then by using induction on $r$, we can prove that $\Tl{H^{2r}}{0}e^{\lambda}=0$ for any $\lambda\in\h$ and $r\geq2$. 
\end{remark}
 
We shall prove that $\Tl{H^{2r}}{0}$ acts diagonally on $\Fremo{\lambda}$ for all $\lambda\in\h$ and $\Fretw{}$:
\begin{proposition}\label{cccc}
Let $r\in\Z_{>0}$. 
 
(1) For any $\lambda\in\h$, $\Tl{H^{2r}}{0}$ acts on $\Fremo{\lambda}$ diagonally, and the eigenvalues for $\Tl{H^{2r}}{0}$ are nonnegative integers. 
Moreover, the eigenspace with eigenvalue zero is $\C e^{\lambda}$.

(2) $\Tl{H^{2r}}{0}$ acts on $\Fretw{}$ diagonally. 
The eigenvalues are rational numbers of the form $q_r+\frac{n}{2^{2r-1}}$ with nonnegative integer $n$ , where $q_r$ is the eigenvalue of $\vactw$ for $\Tl{H^{2r}}{0}$.
\end{proposition} 
\begin{proof}
Let $\lambda\in\h$. 
Then $\Tl{H^{2r}}{0}e^{\lambda}=0$ by Remark \ref{remark2}.
This and Lemma \ref{Lemma17} show that
\[
\Tl{H^4}{0}h(-n_1)\ldots h(-n_r)e^{\lambda}=\left(\sum_{i=1}^{r}n_{i}^{2r-1}\right)h(-n_1)\ldots h(-n_r)e^{\lambda}
\]
for any $n_i\in\Z_{>0}$. 
This proves (1).
We take $q_r\in\C$ so that $\Tl{H^{2r}}{0}\vactw=q_r\vactw$. 
Then the assertion (2) follows from the facts that 
\[
\Tl{H^4}{0}h(-n_1)\ldots h(-n_r)\vactw=\left(q_r+\sum_{i=1}^{r}n_{i}^{2r-1}\right)h(-n_1)\ldots h(-n_r)\vactw
\]
for $n_i\in\frac{1}{2}+\Zpos$.    
\end{proof}
Proposition \ref{cccc} implies that $[\Tl{H^{2r}}{0},\Tl{H^{2s}}{0}]=0$ on every irreducible $\Free{+}$-modules. 
In fact, $\Tl{H^{2r}}{0}(r\in\Z_{>0})$ are mutually commutative on any weak $\Free{+}$-module:
\begin{lemma}\label{Lemma11}
For any $r,s\in\Z_{>0}$, $[\Tl{H^{2r}}{0},\Tl{H^{2s}}{0}]=0$ on any weak $\Free{+}$-module. 
\end{lemma}
\begin{proof}
For any $i\in\Zpos$, $H^{2r}(i)H^{2s}\in\mathcal{S}_{2r+2s-i-1}$. 
By \eqref{SR1} and Lemma \ref{Lemma17}, we see that $H^{2r}(i)H^{2s}$ is a linear combination of vectors $L(-1)^{k}H^{2r+2s-i-k}$ for $0\leq k\leq 2r+2s-i-1$, where we set $H^{2t+1}=0$ for $t\in\Zpos$.
Then the commutativity formula and the $L(-1)$-derivative property show that there exist constants $a_i\in\C$ for $1\leq i\leq r+s-1$ such that 
\begin{align}\label{eqn-3}
[\Tl{H^{2r}}{0},\Tl{H^{2s}}{0}]=\sum_{i=1}^{r+s}a_i\Tl{H^{2i}}{0}
\end{align} 
on any weak $\Free{+}$-module. 
We note that the coefficients $a_i\,(1\leq i\leq r+s)$ are independent of the choice of weak $\Free{+}$-modules.

To determine the coefficients $a_i$ we consider the actions of the both hand side in \eqref{eqn-3} on the vectors $\frac{1}{2}h(-n)^2\1\in\Free{+}$ for $1\leq n\leq r+s$.   
The left hand side acts on these vectors by the scalar zero.
On the other hand $\frac{1}{2}h(-n)^2\1$ is an eigenvector for the action of the right hand side in \eqref{eqn-3} of eigenvalue $\sum_{i=0}^{r+s-1}a_in^{2i-1}$. 
Since the determinant of the matrix $(i^{2j-1})_{1\leq i,j\leq r+s}$ is the nonzero number $(r+s)!\prod_{1\leq i<j\leq r+s}(j^2-i^2)$, we see that $a_i=0$ for any $i$.  
This proves the lemma.
\end{proof}

The explicit forms of $H^4$ and $H^6$ are given by   
\begin{align}
H^4&=\frac{1}{3}h(-3)h(-1)\1-\frac{1}{3}h(-2)^2\1,\label{kkkk1}\\
H^6&=\frac{1}{5}h(-5)h(-1)\1-\frac{13}{10}h(-4)h(-2)\1+\frac{11}{10}h(-3)^2\1.\label{kkkk2}
\end{align} 
The actions of $\Tl{H^2}{0}(=L(0)),\,\Tl{H^4}{0}$ and $\Tl{H^6}{0}$ on $\Omega(M)$ for any irreducible $\Free{+}$-module $M$ are given in the Table \ref{TOP}. 
\begin{table}
\begin{center}
\begin{tabular}{|c|c|c|c|c|c|}\hline
&$\Free{+}$&$\Free{-}$&$\Fremo{\lambda}\,(\lambda\in\h-\{0\})$&$\Fretw{+}$&$\Fretw{-}$\\
\hline
$\Omega(M)$&$\C\1$&$\C h(-1)\1$&$\C e^{\lambda}$&$\C\vactw$&$\C h(-1/2)\vactw$\\
\hline
$\w$&$0$&$1$&$(\lambda,\lambda)/2$&$1/16$&$9/16$\\ 
$H^4$&$0$&$1$&$0$&$-1/128$&$15/128$\\

$H^6$&$0$&$1$&$0$&$1/256$&$9/256$\\ \hline
\end{tabular}
\end{center}
\caption{The action of $L(0),\,\Tl{H^4}{0}$ and $\Tl{H^6}{0}$ on the top levels}\label{TOP}
\end{table}

The explicit description of the commutation relations among $L(m),\,\Tl{H^4}{n}$ and $\Tl{H^6}{l}$ for $m,n,l\in\Z$ are given in Section \ref{appendix}. 
We use the following commutation relations: 
\begin{align}
&[\Tl{H^4}{0},L(m)]=-3m\Tl{H^4}{m}-\frac{m^2(m+1)}{2}L(m),\label{tttt1}\\
\begin{split}\label{tttt2}
&[\Tl{H^4}{0},\Tl{H^4}{m}]\\
&\quad=-3m\Tl{H^6}{m}-\frac{m^2(3m+7)}{4}\Tl{H^4}{m}-\frac{m^2(m-1)(m+2)(m+3)}{180}L(m),\end{split}\\
\begin{split}\label{tttt3}
&[\Tl{H^6}{0},L(m)]\\
&\quad=-5m\Tl{H^6}{m}-\frac{15m^2(m+1)}{4}\Tl{H^4}{m}-\frac{m^2(m+1)(2m^2+2m-1)}{6}L(m)
\end{split}
\intertext{and}
&[\Tl{H^4}{1},L(m)]=-(3m-1)\Tl{H^4}{2}-\frac{m(m+1)(3m+1)}{6}L(m+1)\label{tttt4}
\end{align}
for any $m\in\Z$. 
As one of consequences of these commutation relations, we have the following lemma:
\begin{lemma}\label{yyyy}
Let $M$ be a weak $\Free{+}$-module and $u\in M$ such that $\Tl{H^4}{0}u=hu$ and $\Tl{H^6}{0}u=ku$ for some $h,k\in\C$.
Then 
\[
5(\Tl{H^4}{0}-h)(\Tl{H^4}{0}-h-1)L(1)u+9(\Tl{H^6}{0}-k)L(1)u-L(1)u=0.
\]
\end{lemma}
\begin{proof} 
By Commutation relations \eqref{tttt1} and \eqref{tttt3} with $m=1$, we have 
\begin{align*}
\Tl{H^4}{1}u&=-\frac{1}{3}(\Tl{H^4}{0}-h+1)L(1)u\quad\hbox{and}\\
\Tl{H^6}{1}u&=\left(-\frac{1}{5}(\Tl{H^6}{0}-k)+\frac{1}{2}(\Tl{H^4}{0}-h)+\frac{3}{10}\right)L(1)u.
\end{align*}
Then Commutation relations \eqref{tttt1} and \eqref{tttt2} with $m=1$ yield 
\begin{align*}
(\Tl{H^4}{0}-h)^2L(1)u&=[\Tl{H^4}{0},[\Tl{H^4}{0},L(1)]]u\\
&=[\Tl{H^4}{0},-3\Tl{H^4}{1}-L(1)]u\\
&=9\Tl{H^6}{1}u+\frac{21}{2}\Tl{H^4}{1}u+L(1)u\\
&=-\frac{9}{5}(\Tl{H^6}{0}-k)L(1)u+(\Tl{H^4}{0}-h)L(1)u+\frac{1}{5}L(1)u.
\end{align*}
This shows the lemma.
\end{proof}

We remember that $\{\w,\,J\}$ is a set of generators of $\Free{+}$. 
Since 
\begin{align}\label{JHrelation}
J=-9H^4+4L(-2)^2\1-3L(-4)\1,
\end{align} 
we can see that $\{\w, H^4\}$ is also a set of generators of $\Free{+}$. 
Hence by Proposition \ref{aaee}, we get:
\begin{lemma}\label{lemma45}
Let $M$ be a weak $\Free{+}$-module.
Then $u\in\Omega(M)$ if and only if $L(n)u=\Tl{H^4}{n}u=0$ for any $n\in\Z_{>0}$.
\end{lemma}
\begin{remark}
The relation \eqref{JHrelation} implies that $-9 [H^4]=[J]-[\w]+4[\w]*[\w]$. 
Since $[\w]$ and $[J]$ are generators of Zhu's algebra $A(\Free{+})$ of defining relation \eqref{relations1}--\eqref{relations3}, we see that $[\w]$ and $[H^4]$ are also generators of $A(\Free{+})$ with relations 
\begin{align}
&[\w]*[H^4]-[H^4]*[\w]=0,\label{wHrelations1}\\
&([\w]-1)*([\w]-1/16)*([\w]-9/16)*[H^4]=0,\label{wHrelations2}\\
&(70 [H^4]-132[\w]*[\w]+65[\w]-3)*[H^4]=0.\label{wHrelations3}
\end{align}
\end{remark}
\begin{lemma}\label{iiii}
Let $M$ be a weak $\Free{+}$-module and $N$ a sum of irreducible $\Free{+}$-submodules of $M$.
For $u\in M$, if $L(1)u,\,L(2)u\in N$ and the pair $(\Tl{H^4}{0}u,\Tl{H^6}{0}u)$ is one of pairs $(0,0),\,(-1/128 u,1/256 u)$ or $(15/128 u,9/256 u)$ then $u\in\Omega(M)$.
\end{lemma}
\begin{proof}
Since $N$ is a sum of irreducible $\Free{+}$-modules, $\Tl{H^4}{0}$ and $\Tl{H^6}{0}$ acts on $N$ diagonally by Proposition \ref{cccc}. 
Furthermore, by Lemma \ref{Lemma11}, $N$ is decomposed into a direct sum of simultaneous eigenspaces for $\Tl{H^4}{0}$ and $\Tl{H^6}{0}$. 
If $(h,k)\in\C^2$ is a pair of eigenvalues of a simultaneous eigenvector for $\Tl{H^4}{0}$ and $\Tl{H^6}{0}$ in $N$, then by Proposition \ref{cccc} 
\begin{align}\label{pppp}
(h,k)\in\Zpos\times\Zpos\quad\hbox{or}\quad(h,k)\in\left(-\frac{1}{128}+\frac{1}{8}\Zpos\right)\times\left(\frac{1}{256}+\frac{1}{32}\Zpos\right).
\end{align}

We now write $L(1)u=\sum_{i}u_{i}$ such that each $u_i$ is a simultaneous eigenvector for $\Tl{H^4}{0}$ and $\Tl{H^6}{0}$.
Let $(p,q)$ be one of $(0,0),\,(-1/128,1/256)$ or $(15/128,9/256)$ and assume that $(\Tl{H^4}{0}u,\Tl{H^6}{0}u)=(pu,qu)$. 
By Lemma \ref{yyyy} the pair $(h_i,k_i)$ of eigenvalues of $u_i$ for $\Tl{H^4}{0}$ and $\Tl{H^6}{0}$ satisfies the equation $(5(h_i-p)(h_i-p-1)+9(k_i-q)-1)L(1)u=0.$ 

We claim that $5(h-p)(h-p-1)+9(k-q)-1\neq0$ if $(h,k)\in\C^2$ satisfies one of the cases in \eqref{pppp}.
Suppose that  $5(h-p)(h-p-1)+9(k-q)-1=0$. 
Set $a=h-p-1/2$ and $b=k-q-1/4$, then we have $5a^2+9 b=0$. 
We note that $a\equiv0,\pm1/2^7\pmod{1/2^3\Z}$, that $a^2\equiv0\pmod{1/2^6\Z}$ or $1/2^{14}\pmod{1/2^{9}\Z}$ and that $b\equiv0,\pm1/2^8\pmod{1/2^5\Z}$. 
Since $b\in1/2^8\Z$, we have $a^2\in1/2^6\Z$ and hence $b\in1/2^6\Z$.
Thus $b\in1/2^5\Z$ and so is $a^2$. 
This implies that $a\in1/2^2\Z$, and the case happens only when $h=p$.
Therefore we can conclude $a=-1/2$. 
This is contradiction because $b=-5a^2/9\notin1/2^5\Z$. 

Therefore, $L(1)u=0$ holds, and commutation relations \eqref{tttt1} and \eqref{tttt2} with $m=1$ prove that $\Tl{H^4}{1}u=\Tl{H^6}{1}u=0$.
Now we show that $L(2)u=0$. 
By \eqref{tttt1} with $m=2$ we have $\Tl{H^4}{0}L(2)u=-6\Tl{H^4}{2}u-6L(2)u$.
On the other hand, since $L(1)u=\Tl{H^4}{1}u=0$, \eqref{tttt4} with $m=1$ and the identity $[L(1),\Tl{H^4}{1}]u=0$ imply that $\Tl{H^4}{2}u=-\frac{2}{3}L(2)u$.
Thus we get $\Tl{H^4}{0}L(2)u=-2L(2)u$.
Since every eigenvalues for $\Tl{H^4}{0}$ on $N$ are greater than $-2$, we see that $L(2)u=0$.
Consequently $u$ satisfies that $L(n)u=0$ for any $n\geq1$, hence $\Tl{H^4}{n}u=0$ for any $n\geq1$ by \eqref{tttt4}.  
Then Lemma \ref{lemma45} proves $u\in\Omega(M)$.
\end{proof}

\section{Vanishing pairs of $\Free{+}$-modules for $\Ext{}$}\label{split}
In this section and next section, we argue the splitness of short exact sequences $0\rightarrow M^{1}\rightarrow M\rightarrow M^{2}\rightarrow0$ for irreducible $\Free{+}$-modules $M^{1},\,M^{3}$.
The aim of this section is to find pairs $(M^1,M^2)$ of irreducible $\Free{+}$-modules subject to $\Ext{\Free{+}}(M^2,M^1)=0$.
 
We first prove that $\Ext{\Free{+}}(\Free{\pm},\Free{\pm})=0$ respectively.
\begin{proposition}\label{Prop103}
The extension groups $\Ext{\Free{+}}(\Free{\pm},\Free{\pm})$ are zero respectively.
\end{proposition}
\begin{proof}
Let $M$ be a weak $\Free{+}$-module such that $M$ contains $\Free{+}$ as a weak $\Free{+}$-submodule and that $M/\Free{+}$ is isomorphic to $\Free{+}$.
By Lemma \ref{Prop53}, we see that $M=\bigoplus_{n=0}^{\infty}M_{(n)}$.
Since $\Free{+}$ has no homogeneous vector of weight one, we have $M_{(1)}=0$.
Therefore, $L(-1)u=0$ for any $u\in M_{(0)}$ and hence an $\Free{+}$-submodule generated by a nonzero vector in $M_{(0)}$ is isomorphic to $\Free{+}$ (see \cite{L2}). 
This shows that $M\cong\Free{+}\oplus\Free{+}$. 

Next we consider a short exact sequence $0\rightarrow \Free{-}\xrightarrow{\phi} M\xrightarrow{\psi}\Free{-}\rightarrow 0$ for a weak $\Free{+}$-module $M$.
Since the lowest weight of $\Free{-}$ is $1$, Lemma \ref{Prop53} implies $M=\bigoplus_{n=1}^{\infty}M_{(n)}$ and $M_{(1)}$ is an $A(\Free{+})$-module.   
To prove that $M\cong\Free{-}\oplus\Free{-}$, it suffices to show that $M_{(1)}$ is a completely reducible $A(\Free{+})$-module. 
We take $a=-2^8/105(\w-(1/16)\1)*(\w-(9/16)\1)*H^{4}$. 
Then relations \eqref{wHrelations1}--\eqref{wHrelations3} imply that $[a]$ is an idempotent in $A(\Free{+})$ corresponding to the irreducible module $\Omega(\Free{-})=\C h(-1)\1$, i.e., $[a]$ satisfies  $[\w]*[a]=[a],\,[H^4]*[a]=[a]$ and $[a]*[a]=[a]$.
Since $M_{(1)}$ is two dimensional, it is easy to see that $\Tl{a}{0}|M_{(1)}$ is surjective. 
Hence $M_{(1)}=\image (\Tl{a}{0}|M_{(1)})$ is a direct sum of copies of the irreducible $A(\Free{+})$-module $\C h(-1)\1$.
\end{proof}

The following proposition can be proved by using the same way in Proposition \ref{Prop103}:
\begin{proposition}\label{Prop19}
For any $\epsilon_1,\epsilon_2\in\{\pm\}$, $\Ext{\Free{+}}(\Fretw{\epsilon_1},\Fretw{\epsilon_2})=0$. 
\end{proposition}
\begin{proof}
It is clear that $\Ext{\Free{+}}(\Fretw{\pm},\Fretw{\mp})=0$ from Proposition \ref{Prop73} because the difference of the lowest weights of $\Fretw{+}$ and $\Fretw{-}$ is not an integer.  

Let $\epsilon\in\{\pm\}$ and $0\rightarrow \Fretw{\epsilon}\,\xrightarrow{\phi}\, M\,\xrightarrow{\psi}\, \Fretw{\epsilon}\rightarrow0$ be a short exact sequence for a weak $M$-module. 
By Lemma \ref{Prop53}, $M=\bigoplus_{n\in\Zpos}M_{(n+1/16)}$ (resp. $\bigoplus_{n\in\Zpos}M_{(n+9/16)}$) if $\epsilon=+$ (resp. $-$). 
It suffices to prove that $M_{(1/16)}$ (resp. $M_{(9/16)}$) is a direct sum of copies of the irreducible $A(\Free{+})$-module $\C\vactw$ (resp. $\C h(-1/2)\vactw$). 
As in the proof of Proposition \ref{Prop103}, we have only to find the idempotents $[a_{\pm}]$ in $A(\Free{+})$ corresponding to $\Omega(\Fretw{\pm})$ respectively. 
But we see that they are given by the vectors 
\begin{align*}
a_{+}=-\frac{2^{12}}{15}(\w-\1)*\left(\w-\frac{9}{16}\1\right)*H^4\quad\mbox{and}\quad 
a_{-}=-\frac{2^{12}}{105}(\w-\1)*\left(\w-\frac{1}{16}\1\right)*H^4
\end{align*} 
by using relations \eqref{wHrelations1}--\eqref{wHrelations3}. 
\end{proof}

In the next two proposition, we provide pairs $(M,N)$ such that $\Ext{\Free{+}}(N,M)=0$ for irreducible $\Free{+}$-modules $M$ and $N$ when one of $M$ and $N$ is $\Fremo{\lambda}$ for some $\lambda\in\h$. 
\begin{proposition}\label{qqqq}
For $\lambda,\,\mu\in\h$ with $(\lambda,\lambda)\neq(\mu,\mu)$, $\Ext{\Free{+}}(\Fremo{\lambda},\Fremo{\mu})=0$.
\end{proposition}
\begin{proof}
If $(\lambda,\lambda)/2-(\mu,\mu)/2\notin\Z$, then the proposition follows from Proposition \ref{Prop73}.
By Proposition \ref{Prop7}, we may assume that $(\lambda,\lambda)/2-(\mu,\mu)/2\in\Z_{>0}$ (we note that $\Fremo{\lambda}$ is self dual for any $\lambda\in\h$\,). 
Let $M$ be an extension of $\Fremo{\lambda}$ by $\Fremo{\mu}$. 
Let $\psi:M\to\Fremo{\lambda}$ be a canonical projection. 
By Lemma \ref{Prop53}, $M$ is decomposed into a direct sum of generalized eigenspaces for $L(0),\Tl{H^4}{0}$ and $\Tl{H^6}{0}$.  
Since $L(0),\Tl{H^4}{0}$ and $\Tl{H^6}{0}$ mutually commute by Proposition \ref{Lemma11}, $M=\bigoplus_{r,h,k\in\C}M_{(r,h,k)}$ , where we denote by $W_{(r,h,k)}$ the simultaneous generalized eigenspaces of a weak $\Free{+}$-module $W$ of the eigenvalues $(r,h,k)$ for $L(0),\Tl{H^4}{0}$ and $\Tl{H^6}{0}$. 
Then there exists $u\in M_{((\lambda,\lambda)/2,0,0)}$ such that $\psi(u)=e^{\lambda}$. 
We see that $\Tl{H^4}{0}u,\Tl{H^6}{0}u\in M_{((\lambda,\lambda)/2,0,0)}\cap\Fremo{\mu}$.
Since $M_{((\lambda,\lambda)/2,0,0)}\cap\Fremo{\mu}=\Fremo{\mu}_{((\lambda,\lambda)/2,0,0)}=0$ by Proposition \ref{cccc} (1), we have $\Tl{H^4}{0}u=\Tl{H^6}{0}u=0$. 

Meanwhile, we have $L(1)u,\,L(2)u\in\Fremo{\mu}$ because $\psi(u)=e^{\lambda}\in\Omega(\Fremo{\lambda})$.
Therefore, Lemma \ref{iiii} shows that $u\in\Omega(M)$.
This implies that the $\Free{+}$-submodule $N$ generated by $u$ does not contain $e^{\mu}$.
Therefore, $N\cap\Fremo{\mu}=0$ and then $M=N\oplus\Fremo{\mu}$. 
This proves the proposition.
\end{proof}

\begin{proposition}\label{rrww}
For any $\lambda\in\h$, the extension  group $\Ext{\Free{+}}(\Fretw{\pm},\Fremo{\lambda})$ and $\Ext{\Free{+}}(\Fremo{\lambda},\Fretw{\pm})$ are zero. 
\end{proposition}
\begin{proof}
Since $\Fremo{\lambda}$ and $\Fretw{\pm}$ are self dual, by Proposition \ref{Prop7}, it is enough to show that either $\Ext{\Free{+}}(\Fretw{\pm},\Fremo{\lambda})=0$ or $\Ext{\Free{+}}(\Fremo{\lambda},\Fretw{\pm})=0$ for each $\lambda\in\h$. 
We shall prove this only for $\Fretw{+}$. 
The proof for $\Fretw{-}$ can be done by a similar way.

Let $(M^1,M^2)$ be one of the pairs $(\Fremo{\lambda},\Fretw{+})$ or $(\Fretw{+},\Fremo{\lambda})$, and let $0\rightarrow M^{1}\rightarrow M\rightarrow M^{2}\rightarrow0$ a short exact sequence with a weak $\Free{+}$-module $M$. 
Then by Lemma \ref{Prop53} for $a=H^4$ and Proposition \ref{cccc}, we see that $M$ is a direct sum of eigenspaces for $H^4$. 
This proves that $L(0)$ acts on $M$ diagonally.  

We may assume that $(\lambda,\lambda)/2-1/16\in\Z$.
We separate the proof into two cases: one is the case $(\lambda,\lambda)/2\geq1/16$ and the other is the case $(\lambda,\lambda)/2<1/16$.
In the case $(\lambda,\lambda)/2\geq1/16$, we consider the short exact sequence $0\rightarrow \Fretw{+}\rightarrow M\rightarrow \Fremo{\lambda}\rightarrow0$ for a weak $\Free{+}$-module $M$. 
Let $u$ be a nonzero simultaneous eigenvector for $(L(0),\Tl{H^4}{0})$ of eigenvalues $((\lambda,\lambda)/2,0)$. 
Then we see that $u$ satisfies the condition in Lemma \ref{iiii}.
Hence $u\in\Omega(M)$ and $\C u\cong\C e^{\lambda}$ as $A(\Free{+})$-modules.
Since $\Fretw{+}$ is irreducible, the $\Free{+}$-submodule $N$ generated by $u$ has  the intersection zero with the image of $\Fretw{+}$.
Therefore, $N\cong\Fremo{\lambda}$ and $M\cong\Fretw{+}\oplus \Fremo{\lambda}$.

In the case $(\lambda,\lambda)/2>1/16$, we consider a short exact sequence $0\rightarrow \Fremo{\lambda}\xrightarrow{\phi} M\xrightarrow{\psi} \Fretw{+}\rightarrow0$ for a weak $\Free{+}$-module $M$. 
Let $u$ be a nonzero simultaneous eigenvector for $(L(0),\Tl{H^4}{0})$ of eigenvalues $(1/16,-1/128)$. 
Then $\Tl{H^6}{0}u=1/256 u$ and $\psi(u)\in\C\vactw$. 
Therefore, by Lemma \ref{iiii}, we have $u\in\Omega(M)$. 
This shows that the submodule $N$ generated by $u$ does not contain $\phi(e^\lambda)$.
Since $\Fremo{\lambda}$ is irreducible, we have $N\cap\phi(\Fremo{\lambda})=0$, and $M=N\oplus \Fremo{\lambda}$. 
\end{proof}

Finally we have the following theorem.
\begin{theorem}\label{Theo4}
Let $M^1,M^2$ be irreducible $\Free{+}$-modules. 
Then $\Ext{\Free{+}}(M^2,M^1)=0$ if the pair $(M^1,M^2)$ is one of the following ones: 
\begin{align}
&(\Free{\pm},\Free{\pm}),\label{ses1}\\
&(\Free{\pm},\Fremo{\lambda}),\quad(\Fremo{\lambda},\Free{\pm})\quad\hbox{for $\lambda\in\h-\{0\}$},\label{ses11}\\
&(\Free{\pm},\Fretw{\pm}),\quad(\Fretw{\pm},\Free{\pm}),\label{ese2}\\
&(\Free{\pm},\Fretw{\mp}),\quad(\Fretw{\pm},\Free{\mp}),\label{ese3}\\
&(\Fretw{\pm},\,\Fretw{\pm}),\quad(\Fretw{\pm},\,\Fretw{\mp}),\label{ese8}\\
&(\Fremo{\lambda},\,\Fremo{\mu})\quad\hbox{for $\lambda,\mu\in\h-\{0\}$ such that  $(\lambda,\lambda)\neq(\mu,\mu)$},\label{ese5}\\
&(\Fretw{\pm},\Fremo{\lambda}),\quad(\Fremo{\lambda},\Fretw{\pm}).\label{ese6}
\end{align}
\end{theorem}
\begin{proof}
It is proved in Propositions \ref{Prop103}--\ref{rrww} with Proposition \ref{Prop7} that the theorem is true for the pairs in \eqref{ses1} and \eqref{ese8}--\eqref{ese6}.
As for the pairs in \eqref{ese2} and \eqref{ese3} it follows from Proposition \ref{Prop73}. 
By Proposition \ref{Prop7}, it suffices to show the splitness of short exact sequences $0\rightarrow \Free{\pm}\rightarrow M\rightarrow \Fremo{\lambda}\rightarrow0$ for any weak $\Free{+}$-module $M$ and $\lambda\in\h-\{0\}$.
This sequence naturally induces the short exact sequences $0\rightarrow \Free{}\rightarrow M\oplus\Free{\mp}\rightarrow \Fremo{\lambda}\rightarrow0$. 
Propositions \ref{qqqq}, we have $M\oplus\Free{\mp}\cong\Free{}\oplus \Fremo{\lambda}$.
Hence $M\cong\Free{\pm}\oplus \Fremo{\lambda}$.
This proves that $\Ext{\Free{+}}(\Fremo{\lambda},\Free{\pm})=0$. 
\end{proof}

\section{Indecomposable modules for $\Free{}$}\label{indecomp}
In this section we shall find some pairs of irreducible $\Free{+}$-modules such that the corresponding extension groups are nonzero. 
For this purpose, we discuss the indecomposable $\Zpos$-graded weak $\Free{}$-module $\Free{}[t]$ generated from Zhu's algebra $A(\Free{})(\cong\C[t])$.
Then we will see that for any irreducible modules $M^1,M^2$ such that $\Ext{\Free{+}}(M^2,M^1)\neq 0$, a nontrivial extension $M$ of $M^2$ by $M^1$ is given as a quotient submodule of $\Free{}[t]$.

We define an $\hh^{+}$-module structure on $\C[t]$ by the action 
\[
\rho(h\otimes x^{n})f(t)=t\delta_{n,0}f(t)\quad\hbox{for $h\in\h,\,n\geq0$ and $f(t)\in\C[t]$}
\]
and consider the induced $\h$-module 
\[
\Free{}[t]:=U(\hh)\otimes_{U(\hh^{+})}\C[t]\cong\Free{}\otimes \C[t]\quad(\hbox{linearly as $\C[t]$-modules}).
\] 
This module in fact has a structure of a weak $\Free{}$-module by a canonical way (see \cite{FLM}) and the vertex operator associated to a vector of monomial type is defined by \eqref{action34}.
We define $\Free{}[t](n)=\Free{}(n)\otimes \C [t]$ for any $n\in\Zpos$. 
Then $\Free{}[t]=\bigoplus_{n=0}^\infty\Free{}[t](n)$ is a $\Zpos$-graded $\Free{}$-module.
     
It is known that there is an algebra isomorphism $A(\Free{})\xrightarrow{\sim}\C[t]$ by mapping $[h(-1)\1]$ to $t$ (see \cite{FZ}).  
Therefore, any $A(\Free{})$-module can be seen as a $\C[t]$-module. 
For an $A(\Free{})$-module $W$, we set $\Free{}_{W}:=\Free{}[t]\otimes_{\C[t]} W$ which is a weak $\Free{}$-module.
As it is proved in \cite{FLM}, we have $\Omega(\Free{}_{W})=W$.
Let $W_1,W_2$ be $A(\Free{})$-modules and $f\in\Hom_{A(\Free{})}(W_1,W_2)$.
Then the map $\id\otimes f:\Free{}_{W_1}\to \Free{}_{W_2}$ is an $\Free{}$-module homomorphism. 
This gives a functor from the category of $A(\Free{})$-modules to the category of weak $\Free{}$-modules.
Since $\Free{}_W$ is a free $\C[t]$-module, the functor is in fact exact:  
\begin{proposition}\label{Prop15}
Let $0\rightarrow W_1\xrightarrow{\phi} W_2\xrightarrow{\psi} W_3\rightarrow0$ be a short exact sequence of $A(\Free{})$-modules $W_i\,(i=1,2,3)$. 
Then the short sequence $0\rightarrow M_{W_1}\xrightarrow{\id\otimes\phi} M_{W_2}\xrightarrow{\id\otimes\psi} M_{W_3}\rightarrow0$ is exact. 
\end{proposition}
Propositions \ref{Prop15} shows that for any $A(\Free{})$-modules $W_1\subset W_2$, $\Free{}_{W_1}/\Free{}_{W_2}\cong\Free{}_{W_1/W_2}$. 

For any ideal $I$ of $\C[t]$, we set $M(1)[t]/I:=M(1)_{\C[t]/I}$.
Since $\C[t]/(t-c)\, (c\in\C)$ is a one dimensional irreducible $A(\Free{+})$-module isomorphic to $\C e^{ch}$, Proposition \ref{Prop15} implies the following corollary:
\begin{corollary}\label{Prop12}
For any $c\in\C$, $\Fremo{c h}\cong\Free{}[t]/(t-c)$ as an $\Free{}$-module. 
\end{corollary} 
 
We note that the canonical inclusion $\Free{+}\subset\Free{}$ induces an algebra homomorphism $A(\Free{+})\rightarrow A(\Free{})$.
Hence we can regard any $A(\Free{})$-modules as $A(\Free{+})$-modules via the homomorphism.
More precisely, we see that the homomorphism maps the generators $[\w]$ and $[H^4]$ of $A(\Free{+})$ to $t^2/2$ and $0$ in $\C[t]$ respectively. 
\begin{proposition}\label{Lemma59}
For any nonzero $c\in\C$, the $A(\Free{})$-module $\C[t]/((t-c)^2)$ is an indecomposable $A(\Free{+})$-module.
\end{proposition}
\begin{proof}
Set $v_c=\overline{1}$ and $u_c=\overline{(t-c)}$, where $\bar{a}$ denotes $a+((t-c)^2)\in\C[t]/((t-c)^2)$ for any $a\in\C[t]$. 
Then $L(0)v_c=(c^2/2) v_c+cu_c$ and $L(0)u_c=(c^2/2) u_c$.
This proves $\C[t]/((t-c)^2)$ is not completely reducible because $L(0)$ is not diagonalizable if $c\neq0$. 
\end{proof}
In the proof of Proposition \ref{Lemma59} we have a short exact sequence $0\rightarrow\C e^{ch}\rightarrow  \C[t]/((t-c)^2)\rightarrow\C e^{c h}\rightarrow 0$ of $A(\Free{})$-modules. 
Proposition \ref{Prop15} proves that for $c\in\C-\{0\}$, 
\[
\Fremo{ch}\subset\Free{}[t]/((t-c)^{2})\quad\hbox{and}\quad(\Free{}[t]/((t-c)^{2}))/\Fremo{ch}\cong\Fremo{ch}.
\] 
Therefore, by Proposition \ref{Lemma59}, $\Free{}[t]/((t-c)^{2})$ is an indecomposable $\Free{+}$-module:
\begin{proposition}\label{Propindecomp1}
For any $\lambda\in\h-\{0\}$, $\Ext{\Free{+}}(\Fremo{\lambda},\Fremo{\lambda})$ is nonzero.  
\end{proposition}
 
Define the action of $\theta$ on $\C[t]$ by $\theta(f(t))=f(-t)$ and extend to $\Free{}[t]$ by $\theta(u\otimes f(t))=\theta(u)\otimes \theta(f(t))$ for $u\in\Free{}$ and $f(t)\in\C[t]$. 
Then we have $\theta\circ Y(a,z)\circ\theta=Y(\theta(a),z)$ for $a\in\Free{}$ on $\Free{}[t]$.
Thus, $\theta$ gives an $\Free{+}$-isomorphism from $\Free{}[t]$ to itself. 
If $I$ is a $\theta$-stable ideal of $\C[t]$, then $\theta$ acts on $\Free{}[t]_{I}$ and $\Free{}[t]/I$, and their $\pm1$-eigenspaces for $\theta$ are $\Free{+}$-modules. 
We write these eigenspaces $(\Free{}[t]/I)^{\pm}$ respectively.  

We consider the $\C[t]$-module $\C[t]/(t^2)$ which is closed under the action of $\theta$.
Proposition \ref{Prop15} yields the exact sequence  
\begin{align}\label{sesq}
0\rightarrow \Fremo{0}\rightarrow \Free{}[t]/(t^2)\rightarrow \Fremo{0}\rightarrow0.
\end{align}
of $\Free{+}$-modules. 
In this case, $\C[t]/(t^2)=\C v^+\oplus\C v^-\cong \C \1\oplus\C \1$ as $A(\Free{+})$-modules, where we set $v^{+}=v_0$ and $v^{-}=u_0$ in the proof of Proposition \ref{Lemma59}.
Since $\theta (v^{\pm})=\pm v^{\pm}$, we have 
\begin{align}\label{import1}
(\Free{}[t]/(t^2))^{\pm}=\Free{+}\otimes v_{\pm}\oplus\Free{-}\otimes v_{\mp}.
\end{align}
It is important to remark that the direct sum in \eqref{import1} is just as a vector space but not as an $\Free{+}$-module. 
For example, we see that $L(-1)v^+=h(-1)v^{-}$ and $L(1)h(-1)v^{+}=v^{-}$.
Since the $\Free{+}$-module generated $v^{-}$ is an irreducible $\Free{+}$-module, we get unbefitting short exact sequences 
\begin{align*}
&0\rightarrow\Free{-}\rightarrow (\Free{}[t]/(t^2))^{+}\rightarrow \Free{+}\rightarrow0,\\
&0\rightarrow\Free{+}\rightarrow (\Free{}[t]/(t^2))^{-}\rightarrow \Free{-}\rightarrow0.
\end{align*}
\begin{proposition}\label{Propindecomp2}
The extension groups $\Ext{\Free{+}}(\Free{\pm},\Free{\mp})$ are nonzero respectively.  
\end{proposition}
Combining Propositions \ref{Propindecomp1} and \ref{Propindecomp2}, we have the following theorem:
\begin{theorem}\label{Theo5}
The extension groups 
\begin{align*}
\Ext{\Free{+}}(\Free{\pm},\Free{\mp})\quad\hbox{and}\quad\Ext{\Free{+}}(\Fremo{\lambda},\Fremo{\lambda})\quad\hbox{for $\lambda\in\h-\{0\}$}
\end{align*}
are nontrivial.
\end{theorem}
\section{Rationality of $\charge{+}$; rank one case}\label{rationality}
In this section we complete a proof of rationality of $\charge{+}$ when the rank of $L$ is one. 
One of main results is that for irreducible $\charge{+}$-modules $M^1$ and $M^2$, the extension group $\Ext{\charge{+}}(M^2,M^1)$ is trivial.
Throughout this section we assume that $L$ is of rank one. 

Before discussing short exact sequences for $\charge{+}$-modules, we recall some results related to the structure of $\charge{+}$.
\begin{theorem}\label{theorem73} 

(1) {\rm(\cite[Proof of Theorem 5.13]{DN2})}
Zhu's algebra $A(\charge{+})$ is semisimple. 

(2) {\rm(\cite{Gail})}
The vertex operator algebra $\charge{+}$ is $C_2$-cofinite. 
\end{theorem}
Therefore, by Theorem \ref{extension}, to show the rationality of $\charge{+}$, it suffices to prove that $\Ext{\charge{+}}(M^2,M^1)=0$ for irreducible $\charge{+}$-modules $M^1$ and $M^2$.
Now we start proving of the triviality of $\Ext{\charge{+}}(M^2,M^1)$ for irreducible $\charge{+}$-modules $M^i\,(i=1,2)$.  
First we show the following proposition:
\begin{proposition}\label{Prop9}
Let $M^{1}$ and $M^{2}$ be irreducible $\charge{+}$-modules. 
If the difference of the lowest weights of $M^{1}$ and $M^{2}$ is not a nonzero integer, then $\Ext{\charge{+}}(M^{2},M^{1})=0$.
\end{proposition}
\begin{proof}
Let $r_{1}$ and $r_{2}$ be the lowest weight of the irreducible modules $M^{i}\,(i=1,2)$. 
If $r_1-r_2\notin\Z$, $\Ext{\charge{+}}(M^{2},M^{1})=0$ follows from Proposition \ref{Prop73}. 

Suppose that $r_{1}=r_{2}=r$ and consider a short exact sequence $0\rightarrow M^1\xrightarrow{\phi} M\xrightarrow{\psi} M^2\rightarrow 0$ of $\charge{+}$-modules for a weak $\charge{+}$-module $M$. 
By Lemma \ref{Prop53}, we see that $M=\bigoplus_{n=0}^{\infty}M_{(r+n)}$. 
Then we have an exact sequence $0\rightarrow \Omega(M^{1})\xrightarrow{\phi|\Omega(M^1)} M_{(r)}\xrightarrow{\psi|M_{(r)}} \Omega(M^{2})\rightarrow 0$ of $A(\charge{+})$-modules.
Since Zhu's algebra $A(\charge{+})$ is semisimple, there exists an $A(\charge{+})$-module homomorphism $\eta_{0}:\Omega(M^{2})\to M_{(r)}$ such that $\psi\circ \eta_{0}=\id$. 
Let $N^{2}$ be a $\charge{+}$-submodule of $M$ generated by $\eta_{0}(\Omega(M^{2}))$. 
Since $\Omega((\Ker\psi)\cap N^2)=\Omega(\phi(M^1))\cap\eta_{0}(\Omega(M^{2}))=0$, we see that $(\Ker\psi)\cap N^2=0$.
Hence the restriction $\psi|N^2$ gives an isomorphism of $\charge{+}$-modules.
Therefore, the exact sequence $0\rightarrow M^1\rightarrow M\rightarrow M^2\rightarrow 0$ splits. 
\end{proof}

By Table \ref{Table5} and Proposition \ref{Prop9} we see that the extension groups for the following pairs are trivial:
\begin{align*}
&(M,M)\quad\hbox{for any irreducible $\charge{+}$-module $M$},\quad(\charlam{\frac{\alpha}{2}}^{\pm},\charlam{\frac{\alpha}{2}}^{\mp}),\\
&(\charge{\pm},\charge{T_i,\pm}),\quad(\charge{T_i,\pm},\charge{\pm}),\quad(\charge{\pm},\charge{T_i,\mp}),\quad(\charge{T_i,\pm},\charge{\mp})\quad\hbox{for $i=1,2$}\\
&(\charlam{\frac{\alpha}{2}}^{\pm},\charge{T_i,\pm}),\quad(\charge{T_i,\pm},\charlam{\frac{\alpha}{2}}^{\pm}),\quad(\charlam{\frac{\alpha}{2}}^{\pm},\charge{T_i,\mp}),\quad(\charge{T_i,\pm},\charlam{\frac{\alpha}{2}}^{\mp})\quad\hbox{for $i=1,2$},\\
&(\charge{T_i,\pm},\charge{T_j,\pm}),\quad(\charge{T_i,\pm},\charge{T_j,\mp})\quad\hbox{for $i,j=1,2$}.
\end{align*}
\begin{table}
\begin{center}
\begin{tabular}{|c|c|c|c|c|c|}\hline
$\charge{+}$&$\charge{-}$&$\charlam{r\alpha/2k}\,(1\leq r\leq k-1)$&$\charlam{\alpha/2}^{\pm}$&$\charge{T_i,+}\,(i=1,2)$&$\charge{T_i,-}\,(i=1,2)$\\
\hline
$0$&$1$&$r^2/4k$&$k/4$&$1/16$&$9/16$ \\ 
\hline
\end{tabular}
\end{center}
\caption{The lowest weights of irreducible $\charge{+}$-modules with $(\alpha,\alpha)=2k$}\label{Table5}
\end{table}

The following proposition is essentially proved in \cite{A3}:
\begin{proposition}\label{Pro37}
The extension groups $\Ext{\charge{+}}(\charge{\pm},\charge{\mp})$ are trivial respectively.
\end{proposition}
\begin{proof}
We consider an exact sequence $0\rightarrow \charge{-}\rightarrow M\rightarrow \charge{+}\rightarrow0$ for a weak $\charge{+}$-module $M$.
By Lemma \ref{Prop53}, we see that $M=\bigoplus_{n=0}^{\infty}M_{(n)}$ and $M_{(0)}\cong\C\1$ as $A(\charge{+})$-modules. 
On the other hand, for any nonzero vector $u\in M_{(0)}$, $L(-1)u\in\Omega(M)\cap M_{(1)}\cong\charge{-}$.
Hence by Proposition 4.5 in \cite{A3}, $L(-1)u=0$ in the case $(\alpha,\alpha)\neq2$. 
In the case $(\alpha,\alpha)=2$, we consider the vector $E=e^{\alpha}+e^{-\alpha}\in(\charge{+})_{1}$. 
Then we have $E*E=4\w$.
On the other hand, $\Tl{E}{0}L(-1)u=E(0)L(-1)u=L(-1)E(0)u=0$.
Therefore, $4L(-1)u=4L(0)L(-1)u=\Tl{(E*E)}{0}L(-1)u=E(0)^2L(-1)u=0$.   
In both cases, we have $L(-1)u=0$. 
This implies that $u$ generates a $\charge{+}$-submodule isomorphic to $\charge{+}$ (see \cite{L2}).
Thus we have $M\cong\charge{+}\oplus\charge{-}$ and $\Ext{\charge{+}}(\charge{+},\charge{-})=0$. 
Proposition \ref{Prop7} proves that $\Ext{\charge{+}}(\charge{-},\charge{+})=0$. 
\end{proof}

Finally we prove that $\Ext{\charge{+}}(M^2,M^1)$ is trivial for any remaining pair $(M^1,M^2)$:
\begin{proposition}\label{Prop11}
Set $(\alpha,\alpha)=2k$. 
Then $\Ext{\charge{+}}(M^2,M^1)=0$ for the following pairs $(M^1,M^2)$:
\begin{align*}
(M^1,M^2)=\,&(\charge{\pm},\charlam{\frac{r}{2k}\alpha}),\quad(\charge{\pm},\charlam{\frac{r}{2k}\alpha})\quad\hbox{for $1\leq r\leq k-1$},\\
&(\charge{\pm},\charlam{\frac{\alpha}{2}}^{\pm}),\quad(\charlam{\frac{\alpha}{2}}^{\pm},\charge{\pm}),\quad(\charge{\pm},\charlam{\frac{\alpha}{2}}^{\mp}),\quad(\charlam{\frac{\alpha}{2}}^{\pm},\charge{\mp}),\\
&(\charlam{\frac{r}{2k}\alpha},\charlam{\frac{s}{2k}\alpha})\quad\hbox{for $1\leq r,s\leq k-1$},\\
&(\charlam{\frac{r\alpha}{2k}},\charge{T_i,\pm}),\quad(\charge{T_i,\pm},\charlam{\frac{r\alpha}{2k}})\quad\hbox{for $1\leq r\leq k-1$ and $i=1,2$}.
\end{align*}
\end{proposition} 
We give a proof of Proposition \ref{Prop11}. 
Let $(M^{1},M^{2})$ be one of pairs described in Proposition \ref{Prop11}.
From \eqref{dec1}--\eqref{dec4}, we see that such a pair $(M^1,M^2)$ has the following properties; 
\begin{list}{}{\setlength{\leftmargin}{3ex}
\setlength{\labelwidth}{25ex} }
\item{(1)} one of $M^{i}$ is isomorphic to $\bigoplus_{\beta\in I}\Fremo{\lambda+\beta}$ as an $\Free{+}$-module for some $\lambda\in L^\circ-L$ and subset $I\subset L$, 
\item{(2)} both $M^{1}$ and $M^{2}$ have no common irreducible component as $\Free{+}$-modules. 
\end{list}
Property (1) follows from the fact that one of $M^1$ and $M^2$ is isomorphic to $\charlam{\lambda}$ for $\lambda\in L^{\circ}$ or $\charlam{{\alpha}/{2}}^{\pm}$.
Hence the index set $I$ can be taken to be $L$ or $\Zpos\alpha$. 
Property (2) is clear from irreducible decompositions \eqref{dec1}--\eqref{dec4}.
   
Now we consider a short exact sequence 
\begin{align}\label{ses6}
0\rightarrow M^{1}\,\xrightarrow{\phi}\, M\,\xrightarrow{\psi}\,M^2\rightarrow0
\end{align}
for a weak $\charge{+}$-module $M$. 
By property (1) and Proposition \ref{Prop7}, we may assume that 
\begin{align}\label{eqn5}
M^2\cong\bigoplus_{\beta\in I}\Fremo{\lambda+\beta}
\end{align}
as $\Free{+}$-modules for some $\lambda\in L^{\circ}$. 

Let $N$ be an arbitrary $\Free{+}$-irreducible component of $M^1$.
Then property (2) implies that $N$ is not isomorphic to $\Fremo{\lambda+\beta}$ for any $\beta\in I$.
Therefore, by Theorem \ref{Theo4}, we see that $\Ext{\Free{+}}(N,\Fremo{\lambda+\beta})=0$ for any $\beta\in I$.
Thus, $\Ext{\Free{+}}(M^2,M^1)=0$ by Lemma \ref{lemma3}.
This shows that there is a homomorphism $\eta: M^{2}\rightarrow M$ of $\Free{+}$-modules such that $\psi\circ\eta=\id$.

We denote by $\hat{M}^{2}:=\eta(M^2)$ which is an $\Free{+}$-submodule of $M$ isomorphic to $M^{2}$.
Since $\hat{M}^{2}$ has no intersection with $M^1$ because of property (2), we have $M=\phi(M^1)\oplus \hat{M}^{2}$. 

\begin{lemma}\label{Lemma6}
The $\Free{+}$-module $\hat{M}^{2}$ is closed under the action of $\charge{+}$.
\end{lemma} 
\begin{proof}
We first recall some fusion rules for $\Free{+}$ from \cite{A1}.
For $\beta\in L,\,\mu\in L^{\circ}$ and an irreducible $\Free{+}$-module $N$, the fusion rule of type $\fusion{\Free{+}}{\Fremo{\mu}}{N}$ is nonzero if and only if $N\cong\Fremo{\mu}$, that of type $\fusion{\Fremo{\beta}}{\Fremo{\mu}}{N}$ is nonzero if and only if $N$ is isomorphic to $\Fremo{\mu+\beta}$ or $\Fremo{\mu-\beta}$. 

Let $N^{1}$ and $N^{2}$ be irreducible $\Free{+}$-submodules of $M^1$ and $\hat{M}^{2}$ respectively. 
By property (2), $N^1$ is not isomorphic to $\Fremo{\lambda+\beta}$ for any $\beta\in I$.
We remark that $I$ can be taken to be $L$ or $\Zpos\alpha$. 
The case $I=\Zpos\alpha$ occurs when $M^2\cong\bigoplus_{\beta\in I}\Fremo{\alpha/2+\beta}$.
Then we can prove that for any $\beta\in I$, there is $\beta'\in I$ such that $\Fremo{\lambda-\beta}\cong\Fremo{\lambda+\beta'}$.
Thus for any $\beta\in I$, $N^1$ is not isomorphic to both $\Fremo{\lambda\pm\beta}$.  
This implies that the fusion rule of type $\fusion{\Free{+}}{N^{2}}{N^{1}}$ and of type $\fusion{\Fremo{\gamma}}{N^{2}}{N^{1}}$ is zero for any $\gamma\in L$.
Since $N^1$ and $N^2$ are arbitrary, we see that the fusion rules of type $\fusion{\charge{+}}{\hat{M}^{2}}{M^{1}}$ for $\Free{+}$ is zero. 

Let $Y(\,\cdot\,,z)$ be a module structure of $M$.
Then the restriction of $Y(\,\cdot\,,z)$ gives an intertwining operator of type $\fusion{\charge{+}}{\hat{M}^2}{M^{1}}$ for $\Free{+}$.
More precisely, let $\iota_{2}:\hat{M}^{2}\to M$ be the canonical injection and $p_1:M\to M^1$ be the canonical projection, then the operator $\mathcal{Y}(\,\cdot\,,z)$ defined by $\mathcal{Y}(u,z)v=p_1(Y(u,z)\iota_2(v))$ for $u\in \charge{+}$ and $v\in \hat{M}^{2}$ gives an intertwining operator of indicated type for $\Free{+}$.
Since the fusion rule of type $\fusion{\charge{+}}{\hat{M}^{2}}{M^{1}}$ is zero we have  
$p_1(Y(u,z)\iota_2(v))=0$ for any $u\in \charge{+}$ and $v\in \hat{M}^{2}$.
This shows that all coefficients in $Y(u,z)v$ for $u\in \charge{+}$ and $v\in \hat{M}^{2}$ are in $\hat{M}^2$ because $M=M^1\oplus \hat{M}^2$. 
Therefore, we see that $\hat{M}^2$ is closed under the action of $\charge{+}$.   
\end{proof} 
Consequently, $M$ is isomorphic to a direct sum of $M^1$ and $M^2$ as a $\charge{+}$-module. 
This implies that $\Ext{\charge{+}}(M^2,M^1)=0$, and the proof of Proposition \ref{Prop11} is completed.

By Propositions \ref{Prop9}--\ref{Prop11}, we get the following theorem:
\begin{theorem}\label{splittheorem}
Let $L$ be a rank one positive definite even lattice. 
Then $\Ext{\charge{+}}(M^2,M^1)=0$ for any irreducible $\charge{+}$-module $M^{1}$ and $M^{2}$.
\end{theorem}

Finally we have the rationality of $\charge{+}$ by Theorems \ref{extension}, \ref{theorem73} and \ref{splittheorem}:
\begin{theorem}\label{maintheorem}
For a rank one positive definite even lattice $L$, the vertex operator algebra $\charge{+}$ is rational.
\end{theorem}
\begin{remark}It is know that $\charge{+}$ is isomorphic to $V_{\Z\beta}$ with $(\beta,\beta)=8$ when $L$ is the root lattice od type $A_1$ (\cite{DG}) and isomorphic to $L(1/2,0)\otimes L(1/2,0)$ when $L=\Z\alpha$ with $(\alpha,\alpha)=4$ (\cite{DGH}), where $L(1/2,0)$ is the Virasoro vertex operator algebra of central charge $1/2$. 
In \cite{A3}, the author proved the rationality of $V_{\Z\alpha}^+$ when $(\alpha,\alpha)/2$ is a prime integer.   
\end{remark} 

\section{Rationality of $\charge{+}$; general cases}\label{Generalrank}
In this section we prove the rationality of $\charge{+}$ in the case the lattice $L$ has a general rank.
The key lemma is the following theorem due to Miyamoto (see \cite[Theorem 6.11]{Miya3}).
\begin{theorem}\label{key-theorem}
Let $V$ be a simple vertex operator algebra and $G$ be a finite group of automorphisms. 
If $V^{G}$ is $C_2$-cofinite and rational, then for any $g\in G$ the vertex operator algebra $V^{\langle g\rangle}$ is also $C_2$-cofinite and rational. 
\end{theorem} 

We start the proof the rationality of $\charge{+}$. 
Let $d$ be the rank of $L$ and take a subset $\{\alpha_1,\ldots,\alpha_d\}$ from $L$ so that $(\alpha_i,\alpha_j)=0$ if $i\neq j$. 
Then we have a sublattice $\oplus_{i=1}^{d}\Z\alpha_i\subset L$. 
This induces an embedding 
\[
\otimes_{i=1}^{d}V_{\Z\alpha_i}=V_{\oplus_{i=1}^{d}\Z\alpha_i}\subset\charge{}\]   
of vertex operator algebras. 
We next set $k_i=(\alpha_i,\alpha_i)$ and consider an automorphism $g_i:=\exp(\frac{2\pi i}{k_i}\alpha_{i}(0))$ of $\charge{}$ for $1\leq i\leq d$. 
It is clear that the group $K:=\langle g_1,\ldots,g_d\rangle$ is a finite subgroup of $\Aut(\charge{})$ and that 
\[
\otimes_{i=1}^{d}V_{\Z\alpha_i}=\charge{K}.
\] 
Now we consider the group $\tilde{K}=\langle g_1,\ldots,g_d,\theta\rangle$ of automorphisms. 
Since $\theta\circ g_i\circ \theta^{-1}=g_i^{-1}$, $G\cong\langle\theta\rangle\ltimes K$ which is also finite. 
Set 
\begin{align}\label{G-orb}
U=\charge{\tilde{K}}=(\otimes_{i=1}^{d}V_{\Z\alpha_i})^{\langle\theta\rangle}.
\end{align}
Then we see that $U$ is a direct sum of tensor products $\otimes_{i=1}^{d}V_{\Z\alpha_i}^{\varepsilon_i}$ with signs $\varepsilon_i\in\{\pm\}\,(i=1,\ldots,d)$ such that the number of $i$ with $\varepsilon_i=-$ is even. 

Finally, we denote by $\theta_i$ the automorphism $\id\otimes\cdots\otimes\theta\otimes\cdots\otimes\id$ of $\otimes_{i=1}^{d}V_{\Z\alpha_i}$ for each $i$, where $\theta$ acts on the $i$-the component. 
Then we see that each $\theta_i$ induces an automorphism of $U$ because $\theta_i$ commutes with $\theta$ in view of \eqref{G-orb}.
Therefore, the tensor product $\otimes_{i=1}^{d}V_{\Z\alpha_i}^{+}$ is given as the fixed point set of $U$ by the automorphism group $H:=\langle \theta_1,\ldots,\theta_d\rangle$. 

By Theorem \ref{maintheorem}, the vertex operator algebra $V_{\Z\alpha_i}^{+}$ is rational and $C_2$-cofinite for each $i$. 
Therefore, the tensor product $\otimes_{i=1}^{d}V_{\Z\alpha_i}^{+}$ is also rational and $C_2$-cofinite.
Since $\otimes_{i=1}^{d}V_{\Z\alpha_i}^{+}=U^{H}$ with abelian finite group $H$ of automorphisms of $U$, Theorem \ref{key-theorem} shows that $U$ is rational and $C_2$-cofinite; we apply the theorem to $V=U$ for $G=H$ and $g=\id$. 
Consequently, by using Theorem \ref{key-theorem} for $V=\charge{},\,G=\tilde{K},\,g=\theta$ and the identification \eqref{G-orb}, we get the rationality of $\charge{+}$.
\begin{theorem}
Let $L$ be a positive definite even lattice.
Then the vertex operator algebra $\charge{+}$ is rational and $C_2$-cofinite.
\end{theorem}

\section{Appendix}\label{appendix}
In this section we give the commutation relations among $L(m),\,H^4(n)$ and $H^6(l)$ for $m,n,l\in\Z$, which defined in Section \ref{vector}. 
We use the notations in Section \ref{vector}.
Since the calculations are very complicated and routine, we only explain the way to get the commutation relations among $L(m),\,\Tl{H^4}{m}$ and $\Tl{H^6}{n}]$ for $m,n\in\Z$. 

We first express homogeneous vectors of monomial types $h(-p)h(-q)\1$ in $\bigoplus_{r=0}^{7}\mathcal{S}_r$ as a sum of vectors of the form $\1,\,L(-1)^i\w\,(0\leq i\leq5),\,L(-1)^iH^4\,(0\leq i\leq3)$ and $L(-1)^iH^6\,(0\leq i\leq1)$.
We next write $\w(i)H^4,\,H^4(i)H^4$ and $\w(i)H^6$ for $i\geq0$ as linear combinations of vectors of monomial types.
Since $\w(i)H^4,\,H^4(i)H^4$ and $\w(i)H^6$ for $i\geq0$ are in $\bigoplus_{r=0}^{7}\mathcal{S}_r$, these vectors can be expressed by using vectors of the form $\1,\,L(-1)^i\w\,(0\leq i\leq5),\,L(-1)^iH^4\,(0\leq i\leq3)$ and $L(-1)^iH^6\,(0\leq i\leq1)$.
Finally by using commutativity formula \eqref{commhom} and $L(-1)$-derivative property, we write the commutation relations as a linear combination of $L(l),\,\Tl{H^{3}}{m}$ and $\Tl{H^6}{n}$.

\begin{align*}
&[L(m),\Tl{H^4}{n}]\\
&\quad=(3m-n)\Tl{H^4}{m+n}+\frac{m(m+1)(3m+n)}{6}L(m+n)-\frac{5}{3}\binom{m+1}{5}\delta_{m,-n}\id,
\end{align*}
\begin{align*}
&[\Tl{H^4}{m},\Tl{H^4}{n}]\\
&\quad=3(m-n)\Tl{H^6}{m+n}+\frac{m-n}{12}(9m^2-2mn+9n^2+21m+21n)\Tl{H^4}{m+n}\\
&\qquad+\frac{m-n}{180}(3 m^4+2m^3n+3m^2n^2+2mn^3+3n^4\\
&\qquad+12m^3+11m^2n+11mn^2+12n^3+3m^2+11mn+3n^2-18m-18n)L(m+n)\\
&\qquad+\frac{5}{3}\binom{m+3}{7}\delta_{m,-n}\id,
\end{align*}
\begin{align*}
&[L(m),\Tl{H^6}{n}]=(5m-n)\Tl{H^6}{m+n}+\frac{3m(m+1)(5m+n)}{4}\Tl{H^4}{m+n}\\
&\qquad+\frac{m(m+1)}{120}(40m^3+56m^2n+19mn^2+n^3\\
&\qquad+40m^2+3mn-7n^2-20m-12n)L(m+n)\\
&\qquad+\frac{1}{2}\binom{m+1}{7}\delta_{m,-n}\id 
\end{align*}
for $m,n\in\Z$.


\begin{thebibliography}{00000}
\bibitem[Yam]{Gail} 
G.~Yamskulna, $C_2$-cofiniteness of the vertex operator algebra $V_L^+$ when $L$ is a rank one lattice, math.QA/0202056.
\end{thebibliography}
\end{document}